\documentstyle[12pt,epsf]{article}
\newcommand{\lra}{\longrightarrow}
\newcommand{\complex}{{\bf C}}
\newcommand{\proj}{{\bf P}}
\newcommand{\integer}{{\bf Z}}
\newcommand{\rational}{{\bf Q}}
\newcommand{\real}{{\bf R}}
\newcommand{\cC}{{\cal C}}
\newcommand{\cD}{{\cal D}}
\newcommand{\cE}{{\cal E}}
\newcommand{\cF}{{\cal F}}
\newcommand{\cG}{{\cal G}}
\newcommand{\cH}{{\cal H}}
\newcommand{\cI}{{\cal I}}

\newcommand{\cL}{{\cal L}}
\newcommand{\cO}{{\cal O}}

\newcommand{\cR}{{\cal R}}
\newcommand{\cP}{{\cal P}}
\newcommand{\cM}{{\cal M}}

\newtheorem{theorem}{Theorem}[subsection]
\newtheorem{example}[theorem]{Example}
\newtheorem{proposition}[theorem]{Proposition}
\newtheorem{prop}[theorem]{Proposition}
\newtheorem{lemma}[theorem]{Lemma}
\newtheorem{corollary}[theorem]{Corollary}
\newtheorem{remark}[theorem]{Remark}
\newtheorem{definition}[theorem]{Definition}

\newtheorem{conjecture}[theorem]{Conjecture}
\newcommand{\prf}{{\it Proof:} }
\newcommand{\qed}{\hspace*{\fill}\hbox{$\Box$}}

\begin{document}
\title{Mirror  
Symmetry and Toric Degenerations of Partial 
Flag Manifolds}
\author{{\sc  Victor V. Batyrev}\thanks {
Mathematisches Institut, Eberhard-Karls-Universit\"at  T\"ubingen,  
D-72076 T\"ubingen, Germany, 
{\em email address:} batyrev@bastau.mathematik.uni-tuebingen.de},\and
{\sc Ionu\c t Ciocan-Fontanine }\thanks {
Department of Mathematics, Oklahoma State University, Stillwater, OK 74078, USA 
{\em email address:} ciocan@math.okstate.edu},
\and{\sc Bumsig Kim}\thanks{
Department of Mathematics, University of California Davis, Davis, CA 95616, USA 
{\em email address:} bumsig@math.ucdavis.edu},  
\and 
{\sc Duco van Straten}\thanks {FB 17, Mathematik, 
Johannes Gutenberg-Universit\"at Mainz, D-55099 Mainz, Germany, 
{\em email address:} straten@mathematik.uni-mainz.de}}

\date{}
\maketitle

\thispagestyle{empty}

\begin{abstract} 
In this paper we propose and discuss a mirror construction 
for complete intersections in partial flag manifolds 
$F(n_1, \ldots, n_l, n)$. 
This construction includes our previous mirror construction 
for complete intersection in Grassmannians and the mirror construction 
of Givental for complete flag manifolds.
The key idea of our construction is 
a degeneration of $F(n_1, \ldots, n_l, n)$ 
to a certain Gorenstein toric Fano variety   
$P(n_1, \ldots, n_l, n)$  which has been investigated 
by Gonciulea and Lakshmibai. We describe  
a natural small crepant desingularization of $P(n_1, \ldots, n_l, n)$ 
and prove a generalized version of a conjecture of  
Gonciulea and Lakshmibai 
on  the singular locus of  $P(n_1, \ldots, n_l, n)$. 
\end{abstract}

\newpage

\tableofcontents

\newpage

\section{Introduction}

Using combinatorial dualities for reflexive polyhedra and Gorenstein cones 
together with the theory of generalized GKZ-hypergeometric 
functions, one can extend the calculation of the number $n_d$ of rational
curves of degree $d$ on the generic quintic threefold in $\proj^4$
by Candelas, de la Ossa, Green, and Parkes \cite{candelas} 
to the case of Calabi-Yau complete intersections in toric varieties 
\cite{Ba2,BB2,BS,LB}. 

Another class of examples which includes   Calabi-Yau 
quintic $3$-folds are Calabi-Yau complete intersections in homogeneous Fano 
varieties $G/P$ where $G$ is a semisimple Lie group and $P$ is its 
parabolic subgroup. It is a priori not clear how to find an 
appropriate mirror family for these varieties, because  $G/P$ is 
not a toric variety in general. In 
\cite{grass}, we described a mirror construction for complete interesections
in the Grassmannian $G(k,n)$, which turned out to involve a degeneration
of $G(k,n)$ to a certain singular toric Fano variety $P(k,n)$ introduced 
by Sturmfels in \cite{sturmfels}. 

In this paper we consider the extension of our  
methods to the case of complete intersections in 
{\em arbitrary partial flag manifolds} and give complete 
proofs of statements from \cite{grass}. 

It turns out that the 
Pl\"{u}cker embedding of  any such flag manifold $F:=F(n_1, \ldots, n_l, n)$ 
admits a flat degeneration to a Gorenstein toric Fano variety 
$P(n_1, \ldots, n_l,n)$. This deformation has been 
studied recently  by Gonciulea and Lakshmibai 
in \cite{GL0,GL1, GL2}.
The ``mirror-dual'' toric variety ${\bf P}_{\Delta(n_1, \ldots, n_l,n)}$ 
associated with a reflexive polyhedron $\Delta(n_1, \ldots, n_l, n)$ 
has a nice combinatorial description in terms of a
certain graph $\Gamma:= \Gamma(n_1,\ldots, n_l,n)$ that was introduced 
by Givental for the case of the complete flag manifolds \cite{G2}. 
The idea of toric degenerations has been discussed in more general framework 
in \cite{B}.

Using residue formula, we compute  explicitly a   series 
$\Phi_F:=\Phi_F(q_1, \ldots, q_l)$ associated with the graph $\Gamma$
and conjecture that $\Phi_F$ gives a solution to the quantum 
$\cD$-module associated with 
Gromov-Witten classes and quantum cohomology of the partial flag manifold  
$F$. We note that there is no essential difficulty in checking  the conjecture 
in each particular case at hand, because it involves only 
calculations in the small quantum cohomology ring of 
$F$, for which explicit formulas are known \cite{ICF}. 
Applying the ``trick with factorials'' (see \cite{grass}, 
or \S \ref{ci} below) to a Calabi-Yau complete intersection in 
$F$, we obtain $\Phi_F$ as a specialization of the toric GKZ-hypergeometric 
series, from which the instanton numbers (that is, the virtual numbers of 
rational curves on the 
Calabi-Yau) can be computed via the standard procedure 
(see e.g. \cite{BS}). As the validity of this trick was shown recently 
for general homogeneous spaces \cite{K3}, this implies that 
any instanton numbers 
computed via the usual ``mirror symmetry method'' 
are automatically proven to be correct in all cases for which  our conjecture 
on $\Phi_F$ holds. The series $\Phi_F$  of complete flag manifolds 
has been investigated  by Schechtman \cite{Sch}.  
 
The paper is organized as follows. In Section 2 we introduce main 
combinatorial notions used in the definition of a Gorenstein toric Fano 
variety $P(n_1,\dots ,n_{l},n)$ associated with a given  partial 
flag manifold $F(n_1,\dots ,n_{l},n)$. In Section 3 we investigate 
singularities  of  $P(n_1,\dots ,n_{l},n)$ and show that these singularities 
can be smoothed by a flat deformation to the partial flag manifold
$F(n_1,\dots ,n_{l},n)$. As a consequence of our 
results, we prove a generalized version of a conjecture 
of  Gonciulea and Lakshmibai about the  singular locus of  
$P(n_1,\dots ,n_{l},n)$ 
\cite{GL2}.  
In Section 4 we discuss quantum differential systems 
following ideas of Givental \cite{G1,G2,G3}. Finally, in Section 5 we explain 
the mirror construction for Calabi-Yau complete intersections in partial 
flag varieties $F$ and the  computations of the corresponding hypergeometric 
series $\Phi_F$.

\vskip .2truein 
\noindent{\bf Acknowledgement}
\noindent
We would like to thank 
A. Givental, S. Katz, S.-A. Str{\o}mme, and  E. R{\o}dland 
for helpful discussions and 
the Mittag-Leffler Institute for hospitality. 
The second and third named authors
have been supported by Mittag-Leffler Institute postdoctoral fellowships.

\newpage 

\section{Toric varieties associated with partial flag manifolds  
}\label{toric}

In this section we explain how to associate to an arbitrary  partial 
flag manifold 
$F(n_1,\dots ,n_{l},n)$ certain combinatorial objects: 
a graph $\Gamma(n_1,\dots ,n_{l},n)$, a reflexive polytope 
$\Delta(n_1,\dots ,n_{l},n)$
and a Gorenstein toric Fano variety $P(n_1,\dots ,n_{l},n)$. 

\subsection{The graph $\Gamma(n_1,\dots ,n_{l},n)$}
\vskip 10pt
Let $k_1,k_2,\dots ,k_{l+1}$ be a fixed sequence of positive integers. 
We set  $n_0=0$, $n_i:=k_1+\dots +k_{i}$ $( i =1,\ldots, l+1)$, and   
$n:=n_{l+1}$. Denote by $F(n_1,\dots ,n_{l},n)$  the
partial flag manifold  parametrizing sequences  of subspaces
$$0 \subset V_1 \subset V_2 \subset \cdots \subset V_l 
\subset \complex^n, $$
with $ {\rm dim}\, V_i = n_i \;\;\; (i =1, \ldots, l)$. 
Then 
\[ {\rm dim}\, F(n_1,\dots ,n_{l},n)
= \sum_{i=1}^{l}(n_i-n_{i-1})(n-n_i).  \] 
To symplify notations, we shall often write $F$ 
instead of $F(n_1,\dots ,n_{l},n)$, if there is no confusion about the 
numbers $n_1,\dots ,n_{l},n$. By a classical result of Ehresmann 
(\cite{Ehr}), a natural basis 
for the integral cohomology of $F$ is given by 
the Schubert classes. These are Poincar\'{e} dual 
to the fundamental classes of the closed Schubert cells 
${C}_w \subset F$,
parametrized by permutations $w \in S_n$ modulo the subgroup  
$$ W(k_1,\dots ,k_{l+1}) := S_{k_1} \times \cdots \times S_{k_{l+1}} \subset 
S_n.$$ 
In particular, the Picard group of $F$, which is isomorphic to $H^2(F,\integer)$, is generated 
by $l$ divisors $C_{1}, \ldots ,C_{l}$,
corresponding to the simple transpositions $\tau_i\in S_n$
exchanging $n_i$ and
$n_i+1$.

\begin{definition}
{\rm  Denote by $\Lambda:=\Lambda(n_1,\dots ,n_{l},n)$ the standard 
{\bf ladder diagram} consisting of unit squares (the number of unit squares  in 
$\Lambda$ is equal to the dimension of $F$)
corresponding to  the Schubert cell  
of maximal dimension in the flag manifold  $F$. We place the ladder diagram 
$\Lambda$ in the lower left corner of a $n \times n$-square $Q$.  
The lower left corner of $\Lambda$ (or of $Q$) will be  denoted  
by $O_0$. We denote by  
$O_i$ 
$( i \in \{ 1, \ldots, l \} )$ the common vertex 
of the  diagonal squares $Q_i$ of size   $k_i \times k_i$ and $Q_{i+1}$ 
of  size $k_{i+1} \times k_{i+1}$  
(the figure below illustrates the case $l =4$).} 
\end{definition}

\centerline{Figure 1}

\vskip 20pt
\epsfxsize 8cm
\centerline{\epsfbox{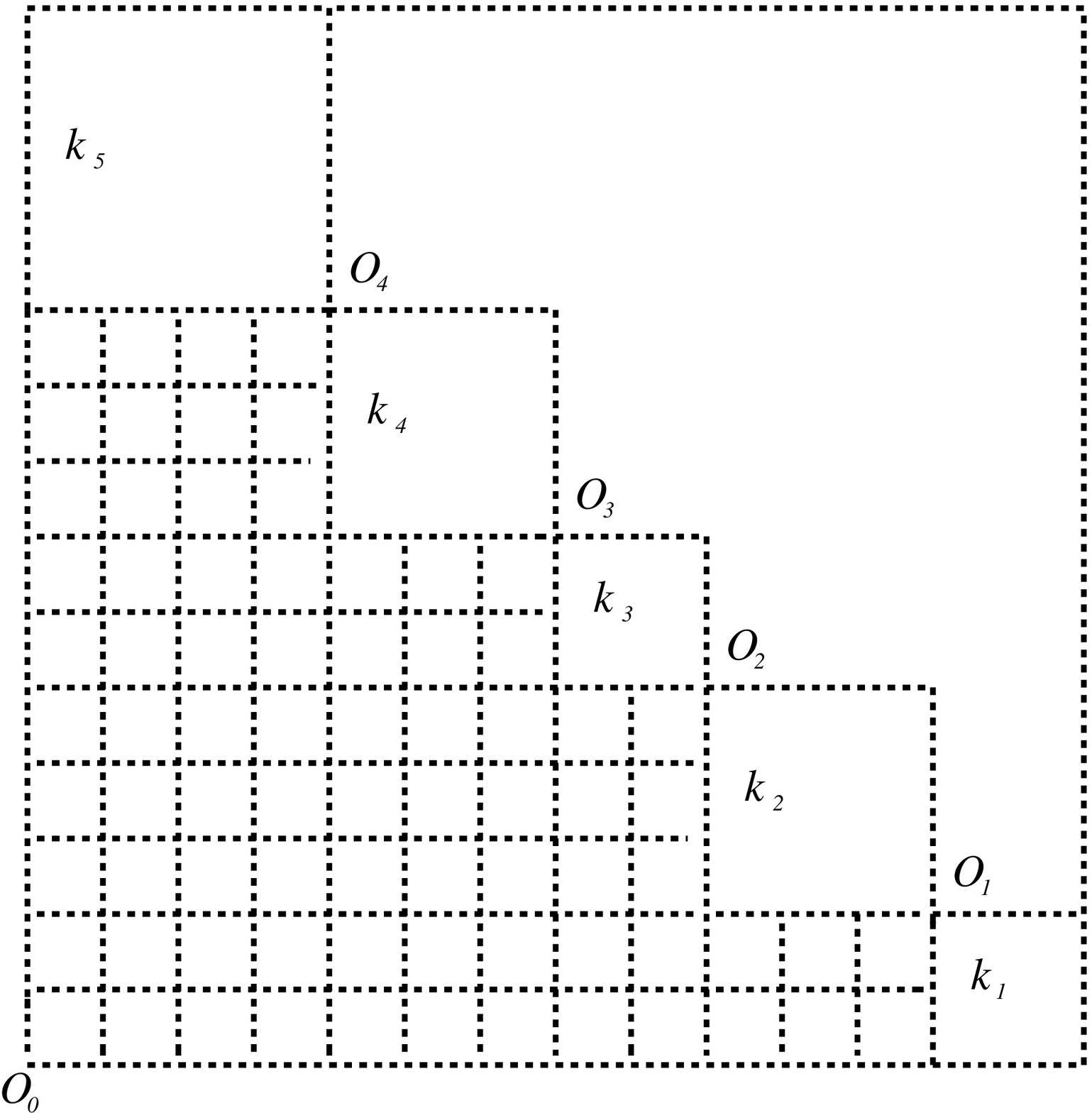}}
\vskip 20pt

\begin{definition}
{\rm Let  $\Lambda = \Lambda(n_1,\dots ,n_{l},n)$ be the above 
ladder  diagram. We {\bf associate} with $\Lambda$ the following: 

\begin{itemize} 
\item $D = D(n_1,\dots ,n_{l},n)$, the set of centers of unit squares  
in $\Lambda$: we place a dot at the center of each unit square and call 
elements of $D$ {\bf dots}.
\item 
$S = S(n_1,\dots ,n_{l},n)$,
the set consisting of $(l+1)$ {\bf  stars}: an element of $S$ 
is obtained by placing a star at the $(1/2,1/2)$-shift 
of the lower left corner of each of the  diagonal squares 
$Q_i$ $(i \in \{1, \ldots, l+1\} )$.
\item $E= E(n_1,\dots ,n_{l},n)$, the set of oriented horizontal and 
vertical segments connecting adjacent elements of $D \cup S$:  
the vertical segments  are oriented downwards, and the horizontal 
segments  are oriented to the right. 
\end{itemize}} 
\end{definition}

\begin{definition}
{\rm $\Gamma :=\Gamma(n_1,\dots ,n_{l},n)$ 
is the oriented graph whose set of vertices is $D \cup S$, 
and whose set of oriented edges is $E$.} 
\end{definition} 

Such a graph $\Gamma$ (without the orientation!) is shown in 
figure 2 below. The edges of $\Gamma$ are drawn
with solid lines. 

\vskip 10pt
\centerline{Figure 2}

\vskip 20pt
\epsfxsize 9cm
\centerline{\epsfbox{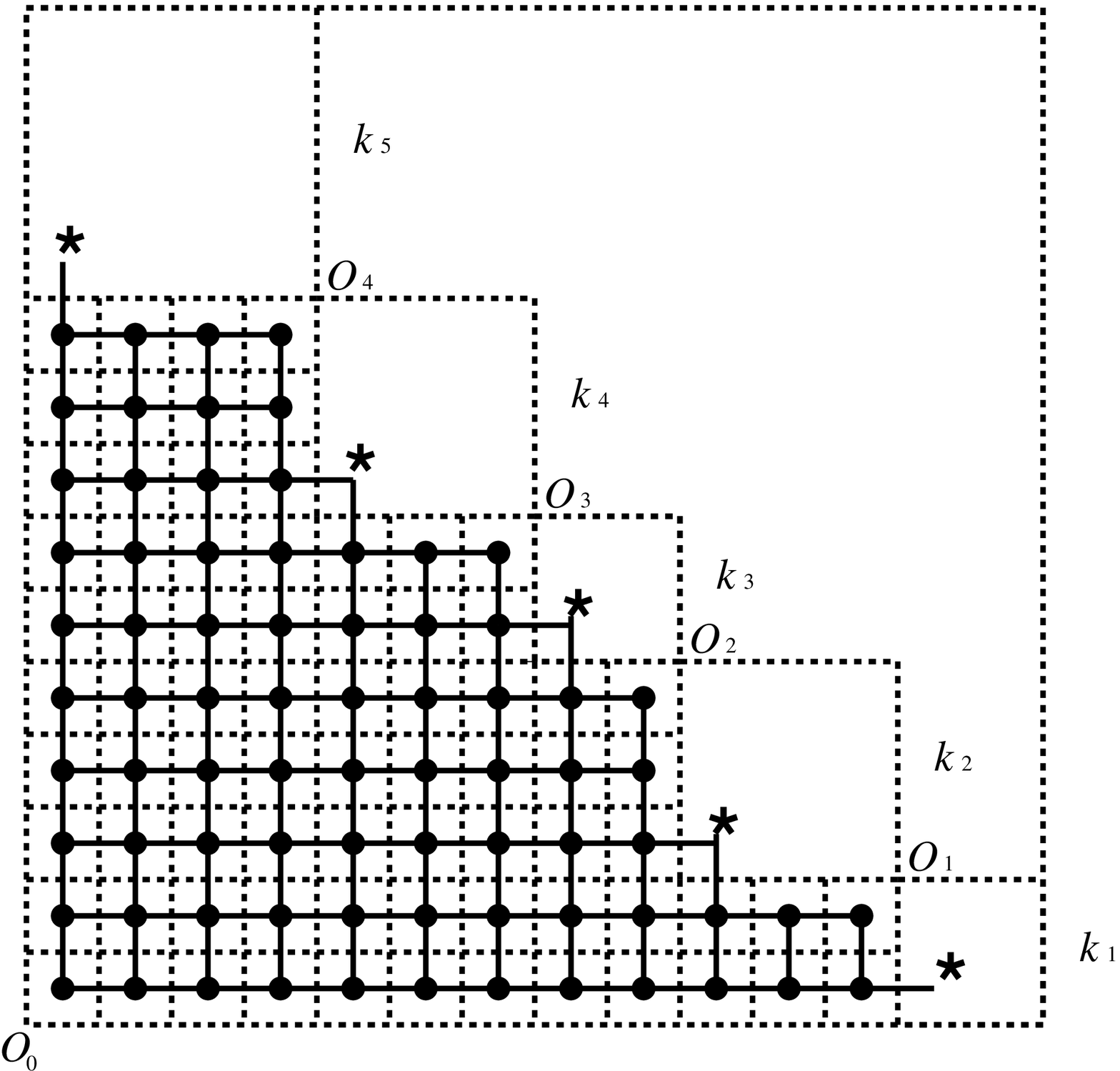}}
\vskip 20pt

\begin{definition}
{\rm We denote by  $L(D) \cong \integer^{\mid D \mid}$, 
$L(S) \cong \integer^{\mid S\mid}$, and  
$L(E) \cong \integer^{\mid E \mid}$ the free abelian groups (or lattices) 
generated by the sets $D, S$ and   $E$.} 
\end{definition}

We remark that the lattices  $L(D) \oplus L(S)$ and $L(E)$ 
can be viewed as the groups of
$0$-chains and $1$-chains of the graph $\Gamma$. 
Then the  boundary map in the chain complex is
$$\partial: L(E) \lra L(D) \oplus L(S)\;\;,\;\;e \mapsto h(e)-t(e),$$
where $h, t:E \lra D \cup S$ are the maps 
that associate to an oriented edge 
$e \in E$ its {\bf head} and its {\bf tail} respectively.

\centerline{Figure 3}

\vskip 20pt
\epsfxsize 4cm
\centerline{\epsfbox{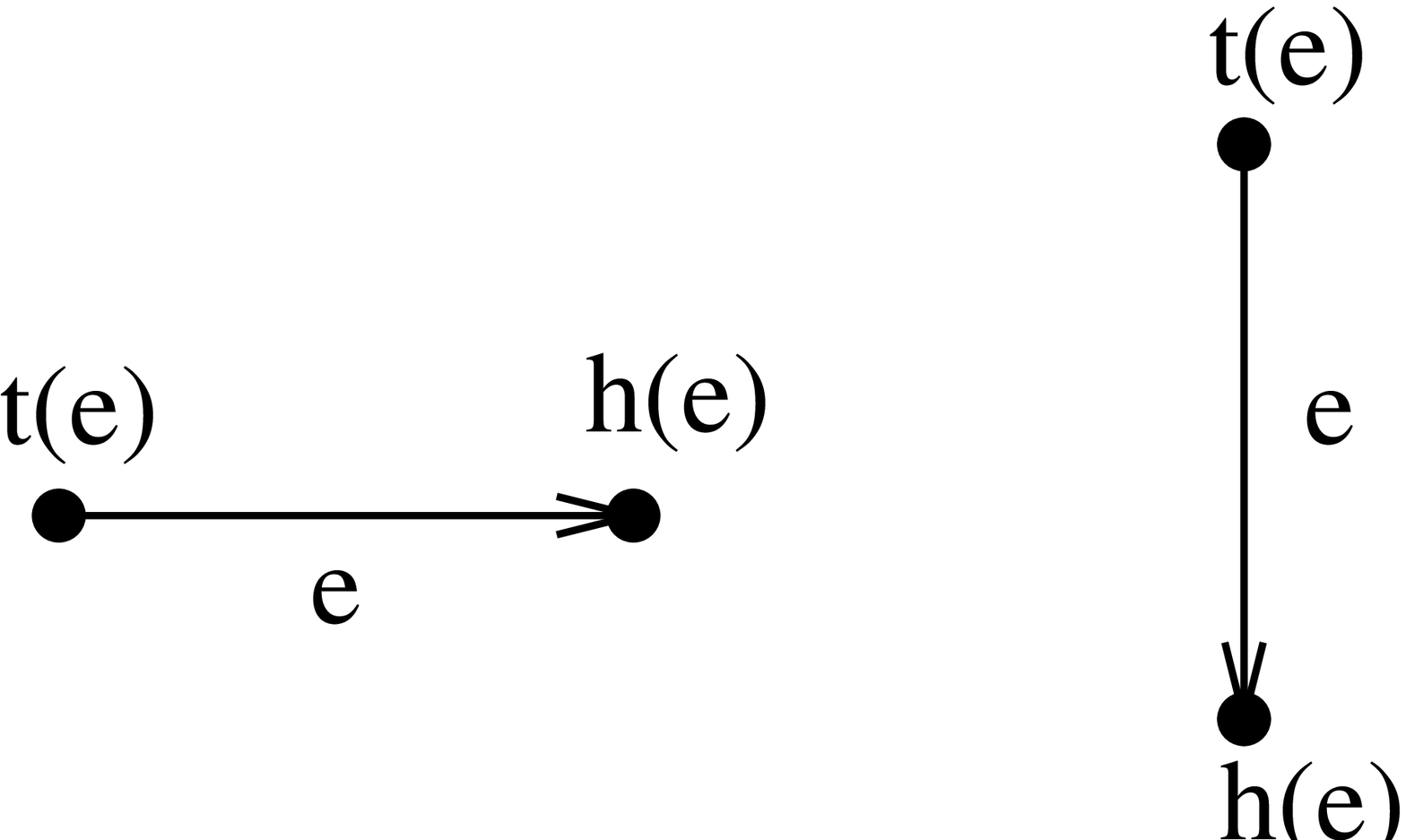}}
\vskip 20pt

\begin{definition}
{\rm A {\bf  box} $b$ in $\Gamma$ is a subset of $4$ edges $\{e,f,g,h \}\subset E$ which form together with their end-points a connected  
subgraph $\Gamma_b \subset \Gamma$, such that the topological space 
associated  to $\Gamma_b$ is 
homeomorphic to a circle. The set of boxes in   $\Gamma$ will be denoted 
by  $B$.}
\end{definition}   

It is easy to see that 

$$
\begin{array}{rccclcrcccl}
H_0(\Gamma)&=&Coker(\partial)&\cong&\integer&\;\;
&H_1(\Gamma)&=&Ker(\partial)&\cong &\integer^{\mid B \mid}
\\
\end{array}
$$   
\vskip 10pt

We also consider the projection $\varrho: L(D) \oplus L(S) \lra L(D)$ and the
composed map $$\delta:=\varrho \circ \partial : L(E) \lra L(D).$$
Since one can regard the groups $L(E)$ and $L(D)$ together with the 
homomorphism $\delta$ as the relative chain complex of the
topological pair $(\Gamma,S)$, we have 

\vskip 10pt

$$
\begin{array}{rccclcrcccl}
H_0(\Gamma,S)&=&Coker(\delta)&=&0&\;\;
&H_1(\Gamma,S)&=&Ker(\delta)&\cong& \integer^{\mid B \mid +l} \\
\end{array}
$$   
\vskip 10pt

\begin{definition}
{\rm  A {\bf  roof}$\,$ $\cR_i$, $i \in \{1,2,\ldots,l\}$ is 
the set of $k_i + k_{i+1}$ edges of $\Gamma$ 
forming the oriented path  
that runs
along the upper right ``boundary'' of $\Gamma$  between  
the $i$-th and the $(i+1)$-st stars in $S$.} 
\end{definition} 

\begin{definition}
{\rm The {\bf  corner} $\cC_b$ of a box  
$b = \{e,f,g,h\} \in B$   
is the pair of edges $\{e,f \} \subset b$ meeting 
at the lower left vertex of $\Gamma_b$. 
So a corner $\cC_b$ contains one vertical edge $e$ 
and one horizontal edge
$f$ such that $h(e)=t(f)$. } 
\end{definition}

The roofs and corners give a decomposition  of the set $E$ of edges of
the graph $\Gamma$ into a disjoint union of subsets: 
$$E=\cR_1 \cup \ldots \cup \cR_l \cup \bigcup_{b \in B} \cC_b.$$
This decomposition is shown in the figure below.

\vskip 10pt
\centerline{Figure 4}

\vskip 20pt
\epsfxsize 8cm
\centerline{\epsfbox{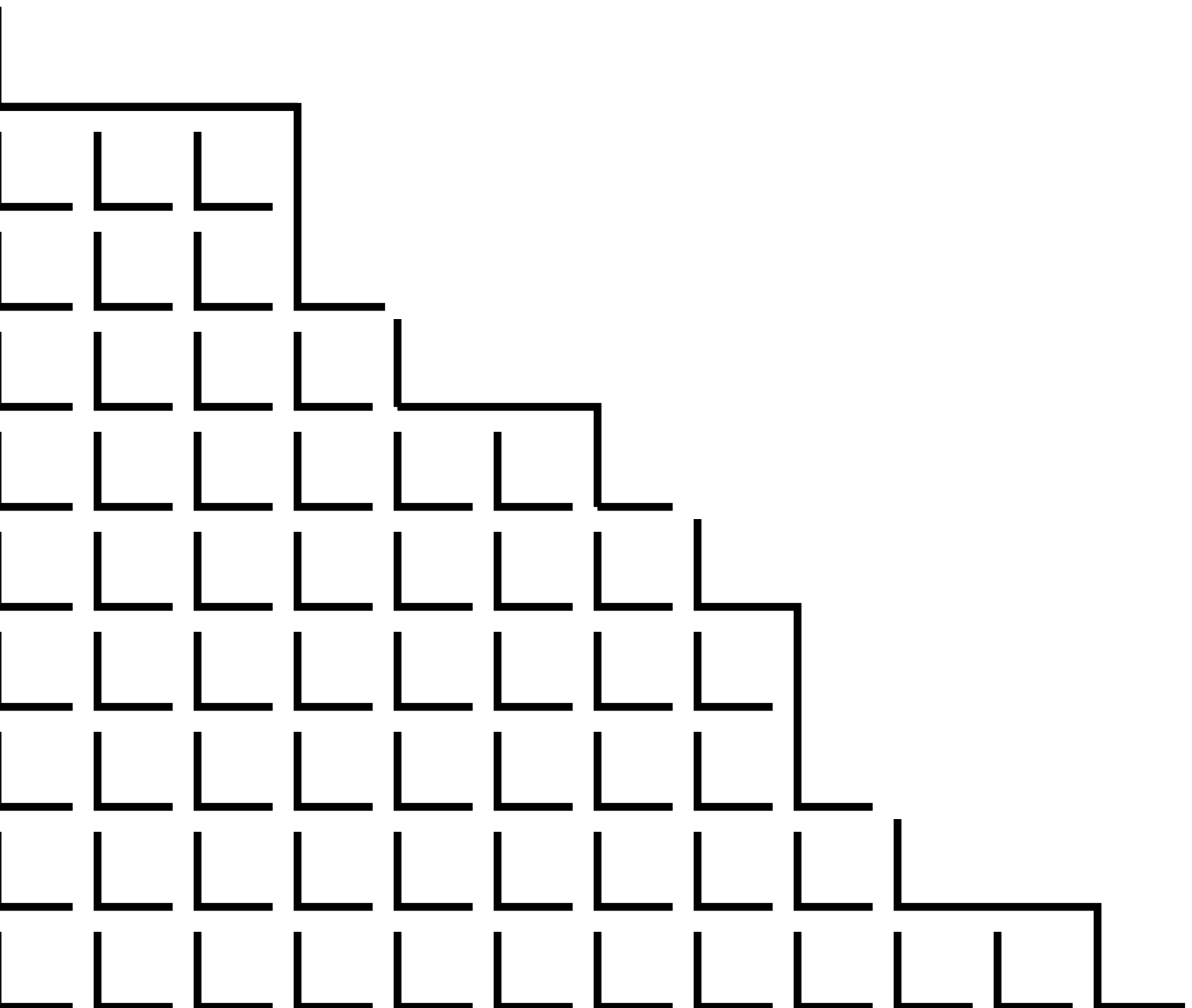}}
\vskip 20pt

\begin{definition}
{\rm The {\bf  opposite corner} $\cC_b^-$ of a 
box  $b = \{ e,f,g,h \} \in B$   
is the pair of edges $\{g,h \} \subset b$ 
meeting at the upper right  vertex of $\Gamma_b$. 
An opposite corner $\cC_b^-$ contains one vertical edge $h$ 
and one horizontal edge
$g$, such that $h(g)=t(h)$. } 
\end{definition}

By elementary arguments one obtains: 

\begin{proposition}
The elements 
\[ \rho_b = \sum_{ e \in \cC_b} e -  \sum_{ e \in \cC_b^-} e, 
\]
where $b$ runs over the set  $B$, 
form a natural ${\bf Z}$-basis of $Ker(\partial) \subset L(E)$.  
Moreover, the elements 
\[ \rho_i = \sum_{e \in \cR_i} e \;\; (i \in \{1, \ldots, l\}) \]
and 
\[ \rho_b = \sum_{ e \in \cC_b} e -  \sum_{ e \in \cC_b^-} e,\;\;
b \in B \]
form a natural ${\bf Z}$-basis of $Ker(\delta) \subset L(E)$.\qed
\label{basis} 
\end{proposition}

\subsection{The toric variety $P(n_1,\dots ,n_{l},n)$}
\vskip 10pt

We denote again by $\delta$ the ${\bf R}$-scalar extension  
$L(E)\otimes\real\lra L(D)\otimes\real$ of the homomorphism 
$\delta\,: \,  L(E) \lra L(D)$.

\begin{definition}{\rm The polyhedron 
$\Delta :=\Delta(n_1,\dots ,n_{l},n)$ associated to $F$ is the convex hull of the set
$$\delta(E) \subset L(D)\otimes\real,$$
where the set $E$ is identified with the standard basis of 
$L(E) \otimes {\bf R} \cong {\bf R}^{\mid E \mid}$.  
}
\end{definition}

In order to describe the faces of the polyhedron $\Delta$ we  introduce
some further combinatorial objects associated to the ladder 
diagram $\Lambda$.

\begin{definition} 
{\rm (i) A {\bf  positive path} $\pi$ in the diagram 
$\Lambda$ is a path obtained by
starting at one of the points $O_i$ $( i =1, \ldots, l)$ 
and moving either downwards, or to the left along
some $n$  edges of $\Lambda$, until the lower left corner $O_0$ 
is reached (see fig. 5).
We denote by $\Pi$ the set of positive paths, and by $\Pi_i$ the set of
positive paths connecting $O_i$ and $O_0$, so that
$$\Pi =\Pi_1 \cup \ldots \cup \Pi_l$$
Note that the number of elements in $\Pi_i$ is $N_i={n \choose n_i}$}

{\rm (ii) A {\bf  meander} is a collection  of positive paths 
$\{ \pi_1, \ldots, \pi_l \}$ $( \pi_i \in \Pi_i)$ , with the property
that the union 
\[ \pi_1 \cup \cdots \cup \pi_l \]
is a {\em tree} with
endpoints $O_0, O_1,\ldots,O_l$. 

The set of all meanders is denoted by
$\cM$.}

\end{definition}

\vskip 10pt
\centerline{Figure 5}

\vskip 20pt
\epsfxsize 8cm
\centerline{\epsfbox{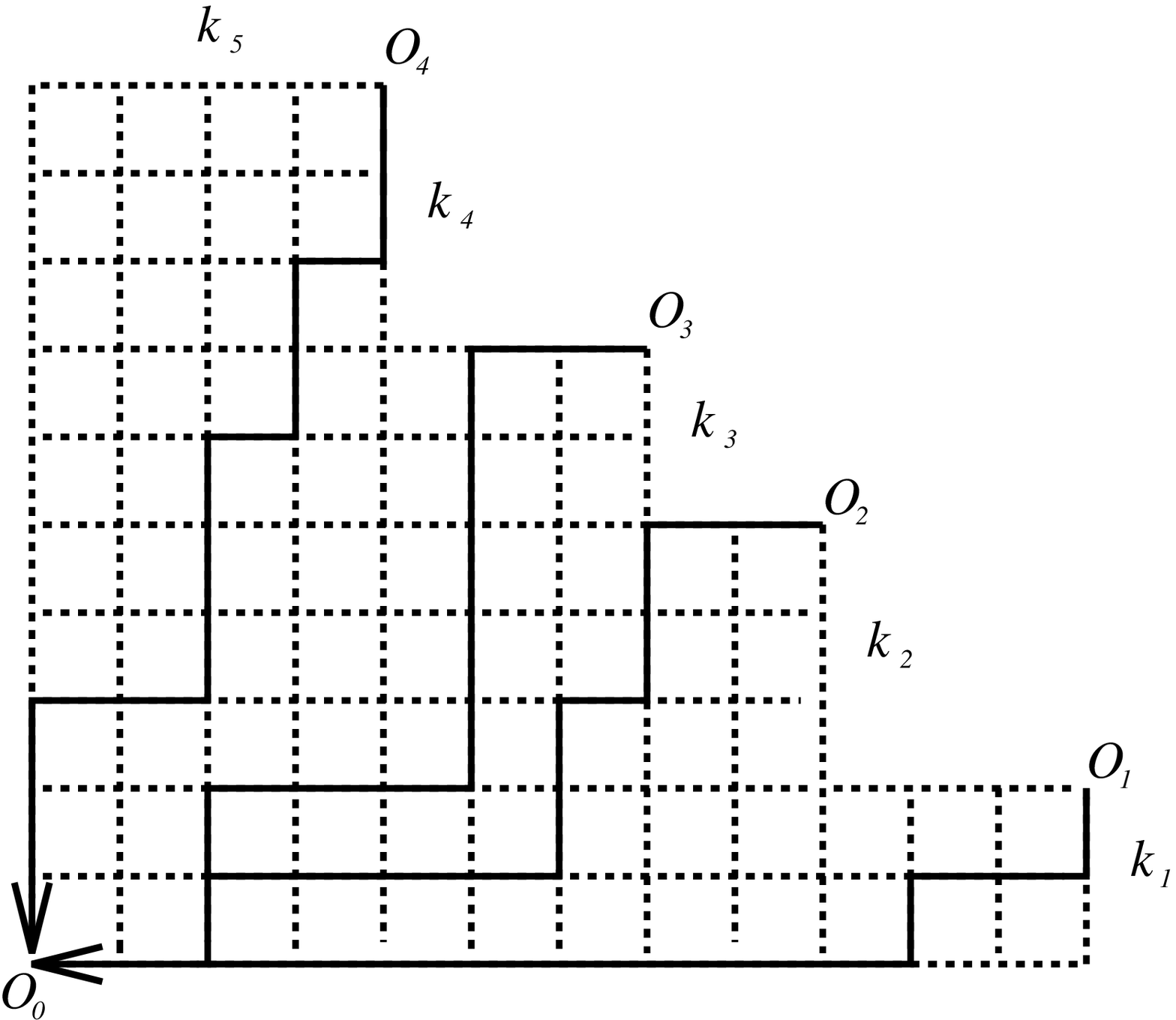}}
\vskip 20pt

\begin{theorem}\label{meande}
There is a natural bijection between the codimension $1$ faces of
$\Delta$ and the set $\cM$ of meanders. 
\end{theorem}
\prf
Since every  face $\Theta$  of $\Delta$ is given by its supporting 
hyperplane,  
it follows from the exact sequence
$$0\lra Ker(\delta)\lra L(E)\otimes\real\stackrel{\delta}\lra 
L(D)\otimes\real\lra 0$$
that this hyperplane can be described by a
linear function 
$$\lambda:L(E)\otimes \real \lra {\real}$$ 
which vanishes on $Ker(\delta)$ and satisfies the conditions  
$$\lambda(v) \le 1,\;\;{\rm for\ all}\; v 
\in L(E)\otimes\real\;\; {\rm with}\;\; \delta(v) \in \Delta ,$$
and 
$$\delta(v) \in \Theta\;\; {\rm iff}\;\; \lambda(v) =  1 \; {\rm and} \;
\delta(v) \in \Delta . $$

Let us show that every meander $m = \{ \pi_1, \ldots, \pi_l \} 
\in \cM$ defines such a linear function $\lambda_m$.
We define the value of $\lambda_m$ on $e \in E$ by the formula: 
\begin{equation}\label{supp}
\lambda_m(e) :=  1-  \sum_{i: \pi_i \cap e \neq \emptyset} 
 \mid \cR_i \mid .
\end{equation} 
It follows that $\lambda_m(e)=1$ if 
the meander $m$ doesn't intersect $e$, and $\lambda_m(e)$ is 
negative if $m$ intersects $e$. Now we show that 
the linear function $\lambda_m$ 
satisfies  the requirement  $\lambda_m|_{Ker(\delta)}=0$.  
By \ref{basis}, it suffices to prove that
\begin{equation}
\sum_{e \in \cR_i} \lambda_m(e)  = 0 , 
\label{r-eqn}
\end{equation}
for all $i \in \{1,2,\ldots,l\}$, and
\begin{equation} 
\sum_{e \in \cC_b } \lambda_m(e) = \sum_{e \in \cC_b^-} \lambda_m(e). 
\label{b-relation}
\end{equation} 
for all $b\in B$.

We remark first 
that every roof $\cR_i$, $i \in \{1,2,\ldots,l\}$ contains exactly 
one edge $e_i \in E$ intersecting the positive path $\pi_i \in m$, 
for which 
$$\lambda_m(e_i) := 1 - \mid \cR_i \mid  < 0.$$ 
On the other hand, $\lambda(e)=1$, for each $e\in\cR_i$, $e\neq e_i$. 
It follows that $$ \sum_{e \in \cR_i} \lambda_m(e)  = 0 \;\;
\forall i \in \{1,2,\ldots,l\}. $$
Now let  $b \in B$ be an arbitrary box. Since the positive 
paths of the meander $m$ form a tree, only the 
following $3$ cases can occur:  

Case 1.  The meander $m$ doesn't intersect edges in $b$. 
Then $\lambda_m(e) =1$ for all $4$ edges of $b$, hence (\ref{b-relation}) holds. 

Case 2. The meander $m$ intersects exactly two edges in $b$. 
Then $m$ intersects exactly one edge $e' \in b$ 
 belonging to $\cC_b$ and exactly 
one edge $e'' \in b$  belonging to $\cC_b^-$. By the formula (\ref{supp}) for $\lambda_m$, we have $\lambda_m(e') = \lambda_m(e'')$. 
So again the relation (\ref{b-relation}) holds. 

Case 3. The  meander $m$ intersects exactly three edges in $b$.
Then $m$ intersects both   edges $e', e'' \in b$ 
 belonging to $\cC_b^-$ and exactly 
one edge $e''' \in b$  belonging to $\cC_b$. By (\ref{supp}) 
$$\lambda_m(e''') = \lambda_m(e') +  \lambda_m(e'') -1.$$ 
Again the relation (\ref{b-relation}) holds. 

Therefore, by Proposition \ref{basis}, $\lambda_m|_{Ker(\delta)}=0$.  

Let $\Theta_m$ be the face of $\Delta$ defined by the supporting 
affine hyperplane $\lambda_m(\cdot) = 1$. We claim that 
$\Theta_m$ has codimension $1$. 
Since $\Theta_m$ is the convex hull of  the  lattice points  
in $\Delta$ 
corresponding to the edges $e\in E$ on 
which $\lambda_m$ takes the
value $1$, it is sufficient to show that any linear function $\lambda'$ 
satisfying $\lambda'|_{Ker(\delta)}=0$ and $\lambda'(v) =1$ for all $v\in L(E)\otimes\real$ with $\delta(v)\in\Theta_m$
must coincide with $\lambda_m$. Indeed, by \ref{basis}, the 
value of such a linear function $\lambda'$ is uniquely determined 
on each edge $e$ of each  roof $\cR_i$ $(1 \leq i \leq l)$: 
$$ \lambda'(e) = \left\{ \begin{array}{ll}  
1 - \mid \cR_i \mid &
\mbox{\rm if $\pi_i \cap e \neq \emptyset$} \\
1 & \mbox{\rm otherwise } \end{array} \right. $$

Next we remark that 
if for some box $b \in B$ we have shown that 
\[ \lambda'(e) = \lambda_m(e) \] 
holds $\forall e \in \cC_b^-$, then, by \ref{basis} and (\ref{b-relation}), 
we obtain 
\[ \sum_{e \in \cC_b } \lambda'(e) = \sum_{e \in \cC_b } \lambda_m(e) \]
and therefore 
\[  \lambda'(e) =  \lambda_m(e) \;\;\forall e \in \cC_b, \] 
since  only one edge $e \in \cC_b$ can be intersected 
by $m$ (see Cases 1-3). Since we have established the equality 
$\lambda'(e) = \lambda_m(e)$ for all $e \in 
\cR_1 \cup \cdots \cup \cR_l$, the 
above arguments imply the equality 
$\lambda'(e) = \lambda_m(e)$ for all $e \in E$.   

Now we prove that any codimension-$1$ face $\Theta$ of $\Delta$ 
can be obtained from some meander $m \in {\cal M}$. For this purpose, 
it suffices to show that if a supporting linear function 
$\lambda$ defines a face $\Theta \subset \Delta$, then there exists 
a meander $m \in {\cal M}$ with $\Theta \subset \Theta_m$. The latter 
is equivalent to the condition $\lambda(e) <1$ for all edges $e \in E$ 
such that $e \cap m \neq \emptyset$.  

First we remark that the linear function $\lambda$ cannot attain the value
$1$ on all edges of the roof $\cR_1$, because $\lambda$ vanishes on the
element $\rho_1 \in  Ker(\delta)$ (see \ref{basis}). 
Now start a positive path $\pi_1$ at $O_1$ whose first nonempty 
intersection with edges of the opposite corner $\cC_{b}^-$ 
of some box $b \in  B$
occurs  on an edge $e_1 \in \cR_1$ with $\lambda(e_1) <1$. 
Since 
\[ \sum_{e \in \cC_b } \lambda(e) = \sum_{e \in \cC_b^-} \lambda(e), \]
the value  of $\lambda$ on at least one of the two 
edges of $\cC_b$  has to be strictly less then $1$. 
We prolong our path through that edge and
enter a next box, where the same reasoning applies. Continuing this, we complete
a positive path $\pi_1$ from $O_1$ to $O_0$ crossing only edges where
$\lambda$ is strictly less than $1$. Now we repeat this construction 
for each of the $O_i$ in
subsequent order, starting at $O_2$ ... etc. If in the process of constructing 
a positive path $\pi_i$ we collide with some already constructed 
positive  path $\pi_j$ $( j < i)$, we just follow from this point 
the path $\pi_j$. 
In the end, we produce a meander with the required  property.  

We conclude that $\Theta_m$ $(m \in {\cal M})$ are all the 
codimension-$1$ faces of $\Delta$.
\qed

\begin{corollary}
$\Delta(n_1,\dots,n_l,n)$ is a reflexive polyhedron.
\label{reflex}
\end{corollary}

\noindent
{\em Proof.} The statement follows immediately from Theorem 
\ref{meande} and from the integrality of the supporting linear function 
$\lambda_m$ (see Definition 4.1.5 in \cite{Ba2}). \qed

\begin{definition}
{\rm The complete rational polyhedral fan $\Sigma = 
\Sigma(n_1,\dots,n_l,n)$ is the fan defined as the collection of cones over 
all faces of $\Delta$. The toric variety 
${\bf P}_{\Sigma}$ associated to the fan $\Sigma$ will 
be denoted by  $P=P(n_1,\dots,n_l,n)$
}
\end{definition}

Using one of the equivalent characterizations of reflexive polyhedra 
(see Theorem 4.1.9 in \cite{Ba2}), we obtain from \ref{reflex}:  

\begin{prop} 
$P(n_1,\dots,n_l,n)$ is a Gorenstein toric Fano variety.\qed
\end{prop}

\section{Further properties of  $P(n_1,\dots,n_l,n)$}

\subsection{Singular locus}
\vskip 10pt

\begin{definition} 
{\rm Define $\widehat{P}=\widehat{P}(n_1,\dots,n_l,n)$ to be the 
toric variety ${\bf P}_{\widehat{\Sigma}}$ associated to the fan 
$\widehat{\Sigma}$, 
obtained by refining the fan $\Sigma$ to a simplicial 
one, whose one-dimensional cones are
the same as the ones of $\Sigma$ (i.e., they are generated by the
lattice vectors $\{\delta(e), e \in E\} \subset L(D)$) and  whose 
combinatorial structure 
is given by the following $|B| +l$ 
primitive collections: 
$$\cR_1, \cR_2, \ldots, \cR_l\; \mbox{\rm and } \;  
\cC_b,\;  b \in B .$$
In other words, the cones of maximal dimension 
of the fan $\widehat{\Sigma}$
are defined by taking all edges $e \in E$  
except  one from each roof and from each corner.} 
\end{definition}

\begin{proposition}
The variety $\widehat{P}$ is a small toric desingularisation
of $P$. 
\label{s-desing}
\end{proposition} 
\prf We have to show that each cone of  $\widehat{\Sigma}$
is contained in a cone of  ${\Sigma}$, and each  cone of 
 $\widehat{\Sigma}$ is generated by a part of a basis. It suffices to prove 
the above properties for cones of $\widehat{\Sigma}$ 
of maximal dimension. 

Choose  an edge $e_i$ in each roof $\cR_i$ $(i =1, \ldots, l)$ and 
an edge $f_b$ in each corner $\cC_b$ $b \in B$. This 
choice determines a $|D|$-dimensional cone $\sigma$ in  $\widehat{\Sigma}$. 
For each 
$i =1, \ldots, l$ there exists a unique positive path from 
$O_i$ to $O_0$ with the following two properties: 
\begin{itemize}
\item  $\pi_i$ crosses the edge $e_i$; 
\item  if $\pi_i$ enters a box $b$, then it  crosses 
the edge $f_b$. 
\end{itemize} 
It is easy to see that the union $\pi_1\cup\ldots\cup\pi_l$ of these paths  
is a meander. Indeed, if a union of positive paths as above is not a tree,
then there must exist a box $b\in B$ with both edges of the corner $\cC_b$
intersecting the union of positive paths. This contradicts the second of
the above conditions.
Therefore the set of edges $\{ e_i \} \cup \{ f_b \}$ 
defines uniquely a meander $m \in {\cal M}$, and the cone 
$\sigma$ is contained in the cone over the face 
$\Theta_m \subset \Delta$. On the other hand, the elements 
$\{ \rho_i \}_{i =1,\ldots, l}$ and $\{ \rho_b \}_{b \in B}$ 
together with the set 
$$G_{\sigma}: = E \setminus \left( \{ e_i \}_{i =1,\ldots, l} \cup 
 \{ f_b \}_{b \in B} \right) $$
form a ${\bf Z}$-basis of $L(E)$. By \ref{basis}, the set 
of generators of $\sigma$ (i.e., 
the $\delta$-image of $G_{\sigma}$)  is a ${\bf Z}$-basis  of 
$L(D)$. 

The desingularization morphism $\widehat{P} \lra P$ induced by 
the refinement $\widehat{\Sigma}$ of $\Sigma$ is small (i.e., 
contracts no divisor), because 
the sets of $1$-dimensional cones in $\widehat{\Sigma}$ and 
$\Sigma$ are the same.    
\qed

\vskip 10pt

There is another way to describe $\widehat{P}$, namely as an 
iterated toric fibration over $\proj^1$:  
One starts with the product of
projective spaces 
\[ {\bf P}^{\mid \cR_1 \mid -1} \times \cdots \times 
 {\bf P}^{\mid \cR_l \mid -1} \]
corresponding to the roofs. Then one chooses  
a corner $\cC_b$ of a box $b \in B$ whose opposite corner $\cC_b^-$ belongs 
to a roof. This choice allows to define a toric bundle over $\proj^1$ 
with the fibre ${\bf P}^{\mid \cR_1 \mid -1} \times \cdots \times 
 {\bf P}^{\mid \cR_l \mid -1}$. Then one adds a  new corner $\cC_{b'}$ 
of a box $b' \in B$ whose 
opposite corner $\cC_{b'}^-$ is contained in the union of roofs and 
$\cC_b$ $\ldots $ etc. 
At each stage of this process one gets a toric fibre bundle over 
$\proj^1$, with 
fibre the space constructed in the previous step. 
Using this description of $\widehat{P}$, one obtains an alternative  
proof  of the fact that
the anti-canonical divisor on $P$ is Cartier and ample, i.e., 
that the polyhedron $\Delta$ is reflexive.

\begin{definition} 
{\rm 
Let $b \in B$ be an arbitrary box. Define $W_b \subset P$ to be the closure 
of the torus orbit in $P$ corresponding to the $3$-dimensional 
cone $\sigma_b$ generated by the $\delta$-image of  the $4$-element set $b$.} 
\end{definition}

\begin{theorem}
The singular locus of $P$ consists of codimension $3$ strata $W_b$, $b \in B$.
These are conifold strata, i.e. transverse to a generic point of $W_b$ the
variety $P$ has an ordinary double point. 
\label{sing}
\end{theorem}
\prf 
Since the desingularization morphism  $\varphi\,:\, 
\widehat{P} \lra P$ is small, 
$P$ is smooth in codimension $2$. Moreover, the singular locus 
of $P$ is precisely the union of toric strata in $P$ over which 
the morphism $\varphi$ is not bijective. According to 
the main result of \cite{Reid}, the exceptional locus $Ex(\varphi) \subset 
\widehat{P}$ (i.e., $\varphi^{-1}(Sing(P))$) is the union 
of toric strata covered by rational curves contracted by $\varphi$. 
On the other hand, since $\widehat{P}$ is an iterated toric bundle, 
the Mori cone $\overline{NE}(\widehat{P})$ is 
a simplicial cone generated by the classes of the  
primitive relations
\[ \sum_{e \in \cR_i } \delta(e) = 0 \;\; i =1, \ldots, l \]
and 
\[ \sum_{e \in \cC_b} \delta(e) = \sum_{e \in \cC_b^-} \delta(e), 
\;\; b \in B \]
(see sections 2 and 4 in \cite{Ba1}).  
Since the morphism $\varphi$ is defined by the 
semi-ample anticanonical class of $\widehat{P}$, it contracts 
exactly the extremal rays in $\overline{NE}(\widehat{P})$ defined 
by the primitive relations corresponding to the boxes $b \in B$. 
The rational curves representing each such class 
cover the  codimension-$2$ strata
$\widehat{W_b},\; b\in B$, corresponding to the 
$2$-dimensional cones 
in $\widehat{\Sigma}$ spanned by the $\delta$-images of the edges forming 
the opposite  corner  $\cC_b^-$.
These strata are contracted, 
with $\proj^1$-fibres, to the codimension $3$ strata $W_b$ in $P$
corresponding to the $3$-dimensional cones $\sigma_b \in \Sigma$ 
over the  quadrilateral faces $\Theta_b$ of $\Delta$ whose vertices 
are $\delta$-images of the edges in $b$ ($b \in B$). It follows that 
$\bigcup_{b \in B}  W_b$ is exactly the singular locus of $P$.  
\qed 

\vskip 10pt

\vskip 10pt

\subsection{Canonical flat smoothing} 

Let $F = F(n_1,\dots,n_l,n)$ be a partial flag manifold. 
The semiample line bundles ${\cal O }(C_{1}), \ldots ,  
{\cal O}(C_l)$ associated to the Schubert divisors $C_1, \ldots, C_l$ 
define the Pl\"ucker embedding 
of $F$ into a product of projective spaces:
$$\phi:F \hookrightarrow {\proj}^{N_1-1}\times \cdots \times 
{\proj}^{N_l-1},$$
where $N_i={n \choose n_i}$. We will always consider  
$F$ as a smooth projective variety together 
with this embedding.

We describe now an embedding of $P$ in the same product of  
projective spaces.

\begin{definition}
{\rm For each $e\in E$, let $H_e$ be the toric Weil divisor on $P$ determined by the $1$-dimensional cone of $\Sigma$ spanned by the vector $\delta(e)$. 

For every edge $e\in \cup_{i=1}^l\cR_i$ which is part of a roof,
denote by $U(e)$ the subset of $E$ consisting of the edge $e$, together with
all edges $f\in E$ which are either directly below $e$ in the graph $\Gamma$, if $e$ is horizontal, or directly to the left of $e$, if $e$ is vertical.}
\end{definition}

Fix $1\leq i\leq l$. For $e\in\cR_i$ consider the Weil divisor $\sum_{f\in U(e)}H_{f}$.  

\begin{lemma}\label{cartier}
For each $e\in\cR_i$, the Weil divisor $\sum_{f\in U(e)}H_{f}$ is Cartier. Moreover, if $e'\in \cR_i$ is another
edge in the same roof, then the associated divisor $\sum_{f'\in U(e')}H_{f'}$ is linearly equivalent to $\sum_{f\in U(e)}H_{f}$.
\end{lemma}

\prf To each edge $e\in \cR_i$, and each positive path $\pi\in\Pi_i$, joining $O_i$ with $O_0$, we associate a linear function $\pi[e]:L(E)\lra \integer$ defined by
$$\pi[e](g)=\left\{ \begin{array}{rl}
 0, &\mbox {\rm if $\pi\cap g=\emptyset$ and $g\notin U(e)$,}\\
 0, &\mbox  {\rm if $\pi\cap g\neq\emptyset$ and $g\in U(e)$,}\\
-1,&\mbox  {\rm if $\pi\cap g=\emptyset$ and $g\in U(e)$,}\\
 1,&\mbox  {\rm if $\pi\cap g\neq\emptyset$ and $g\notin U(e)$,}\end{array} \right. $$

It is an elementary exercise to check that $\pi[e]$
vanishes on the elements
\[ \rho_j = \sum_{g \in \cR_j} g \;\; (j \in \{1, \ldots, l\}) \]
and 
\[ \rho_b = \sum_{ g \in \cC_b} g -  \sum_{ g\in \cC_b^-} g,\;\;
b \in B. \]
It follows from Proposition \ref{basis} that $\pi[e]$ descends to a linear function on $L(D)$.

To show that $\sum_{f\in U(e)}H_{f}$ is Cartier, it suffices to construct for each maximal dimensional cone $\sigma$ in $\Sigma$ an integral linear function $$\lambda_{\sigma}:L(E)\lra \integer ,$$ which vanishes on $ker(\delta)$, and satisfies
\begin{itemize}
\item $\lambda_{\sigma}(g)=0$, for all $g\in E$ such that $\delta(g)\in\sigma$ and $g\notin U(e)$.
\item $\lambda_{\sigma}(g)=-1$, for all $g\in E$ such that $\delta(g)\in\sigma$ 
and  $g\in U(e)$.
\end{itemize}

By Theorem \ref{meande}, every maximal cone $\sigma$ is determined by a meander
$m=(\pi_1,\pi_2,\ldots,\pi_l)$, and $\delta(g)\in\sigma$ if and only if the meander does not intersect the edge $g$ (cf. the proof of \ref{meande}). It follows that $$\lambda_{\sigma}:=\pi_i[e]$$ satisfies the above conditions,
where $\pi_i$ is the positive path in $m$ which joins $O_i$ with $O_0$. Hence $\sum_{f\in U(e)}H_{f}$ is Cartier.
Note that the functional $\lambda_{\sigma}$ defined above does not depend on the
positive paths $\pi_j\; (j\neq i)$ in $m$ that do not intersect the roof $\cR_i$.

To prove the second part of the lemma, define an integral linear function $\mu:L(E)\lra \integer$, by $$\mu(g)=\left\{ \begin{array}{rl}
-1,&\mbox {\rm if $g\in U(e)$,}\\
1, &\mbox {\rm if $g\in U(e')$,}\\
0, &\mbox {\rm otherwise.}\end{array} \right. $$
 
As above, one can easily check that $\mu$ vanishes on $ker(\delta)$, hence it descends to a linear function on $L(D)$. The descended linear function defines a rational function on $P$, whose divisor is $\sum_{f\in U(e)}H_{f}-\sum_{f'\in U(e')}H_{f'}$. This finishes the proof of the Lemma.\qed

\begin{definition} {\rm For each $i=1,2,\ldots,l$, the line bundle associated to 
the roof $\cR_i$ is $$\cL_i:=\cO(\sum_{f\in U(e)}H_{f}),$$ for some edge $e\in\cR_i$.}\label{linebundles}
\end{definition}

It follows from Lemma \ref{cartier} that $\cL_i$ does not depend on the choice of the edge $e\in\cR_i$.

We note that for each maximal dimensional cone $\sigma$ the linear function $\lambda_{\sigma}$ defined in the proof of \ref{cartier} satisfies $\lambda_{\sigma}(g)\geq 0$ for all $g \in E$ such that $g\notin\sigma$. This implies that the line bundle $\cO(\sum_{f\in U(e)}H_{f})$ is generated by global sections (cf. \cite{fulton}, p.68). We will now identify the space of global sections.   

The Cartier divisor $\sum_{f\in U(e)}H_{f}$ determines a rational convex polyhedron $\Delta[e]$ in the dual vector space $L(D)^*\otimes\real$, given by
$$
\Delta[e]=\left\{\lambda\in L(D)^*\otimes\real\; :\; \begin{array}{l}\lambda(\delta(g))\geq -1,\;
\forall g\in U(e),\\ \lambda(\delta(g))\geq 0,\; \forall g\in E\setminus U(e)
\end{array}\right\}.
$$ 
The space of global sections of the line bundle $\cO(\sum_{f\in U(e)}H_{f})$ has a natural basis, indexed by the lattice points in $\Delta[e]$.
By its very definition, for each positive path $\pi\in\Pi_i$, the linear function
$\pi[e]$ introduced in the proof of \ref{cartier} gives such a lattice point. 

\begin{proposition}\label{global-sect}  
For each $i=1,2,\ldots,l$, the space of global sections of $\cL_i$ has a 
natural basis parametrized by the set $\Pi_i$ of positive 
paths connecting $O_i$ and $O_0$.
\end{proposition}

\prf Choose an edge $e\in\cR_i$. We have to show that the only lattice points in $\Delta[e]$ are the ones given by $\pi[e]$, $\pi\in\Pi_i$. Let $\lambda :L(E)\lra\integer$ be any linear function vanishing on $ker(\delta)$, and such that the descended linear function is in $\Delta[e]$. 
 
Since on the one hand $\lambda$ vanishes on every
$$ \rho_j = \sum_{g \in \cR_j} g \;\; (j \in \{1, \ldots, l\}), $$ and on the other hand $\lambda$ can be negative only on edges in $U(e)$, there are exactly two possibilities:

$({\rm I})$ $\lambda(g)=0$ for all $g\in\cup_{j=1}^l\cR_j$;

$({\rm II})$ $\lambda(e)=-1$, there exists an edge $h$ in the $i^{\rm th}$ roof $\cR_i$ with $\lambda(h)=1$, and $\lambda(g)=0$ for all $g\in\cup_{j\neq i}\cR_j\setminus\{e,f\}$. 

If $({\rm I})$ holds, then we start a positive path $\pi$ at $O_i$ that intersects the roof $\cR_i$ at the edge $e$. Let $b$ be the box containing $e$ in its opposite corner and let $f$ be the other edge in $U(e)$ contained in this box. If
$\lambda(f)=0$, we prolong the path through the edge $f$, and enter a next box
$b'$, where we have the same situation as before (i.e., there is an other
edge $f'\in U(e)$, and if $\lambda(f')=0$, then we prolong the path through $f'$ etc.). So we may assume that $\lambda(f)=-1$. The edge $f$ is part of the corner $\cC_b$ of $b$. Let $f''$ be the other edge in $\cC_b$. Since $\lambda$ vanishes on all elements
$$\rho_b = \sum_{ g \in \cC_b} g -  \sum_{ g\in \cC_b^-} g,\;\;b \in B,$$
$\lambda(f'')$ must be strictly positive (hence at least $1$). We prolong the
path $\pi$ through the edge $f''$, and enter a next box $b''$, for which $f''$ is part of the opposite corner. Now $\lambda$ is nonnegative on all four edges of $b''$, and $\lambda(f'')\geq 1$. It follows that there must be an edge
$f'''$ in the corner $\cC_{b''}$, with $\lambda(f''')\geq 1$. We prolong the path through this edge, and enter a next box, where the same reasoning applies. 
Continuing this, we complete eventually a positive path $\pi$. Consider the linear function $\nu:=\lambda-\pi[e]$ on $L(E)$. By construction, and the definition of $\pi[e]$, the functional $\nu$ is nonnegative on all edges $g\in E$. On the other hand, $\nu$ vanishes on the generators of $ker(\delta)$ described in \ref{basis}, since both $\lambda$ and $\pi[e]$ do.
We claim that $\nu$ is identically zero on $L(E)\otimes\real$. Indeed, since $\nu$ is nonnegative and $\nu(\sum_{g\in\cR_j}g)=0$ for $j=1,2,\ldots,l$, it
follows that $\nu$ takes the value zero on each edge in the union of all roofs.
Similarly, if $\nu$ vanishes on each of the edges of the opposite corner $\cC_b^-$ of some box, then it must vanish on each of the edges of the corner $C_b$ as well. From these two facts, one obtains inductively that $\nu$ takes the value zero on every $g\in E$. Hence $\lambda=\pi[e]$.

Assume now that $({\rm II})$ holds. In this case, we start a positive path $\pi$ at $O_i$ that intersects the roof $\cR_i$ at the edge $h$. A reasoning entirely similar to that in case $({\rm I})$ shows that the path $\pi$ can be completed such that the functional $\lambda-\pi[e]$ is nonnegative on every edge. Hence
we obtain again $\lambda=\pi[e]$.\qed

\vskip 10pt
\begin{definition}  
{\rm The $(|D |  +l)$-dimensional 
cone $C = C(n_1,\dots,n_l,n)$ associated to the 
flag manifold $F$ is the convex polyhedral cone 
in the space $Im(\partial)\otimes\real$ spanned by the vectors
$$ \partial(e) \in Im(\partial)\otimes\real 
\cong \real^{\mid D\mid+\mid S\mid-1}$$
with $e \in E$.  We denote by $C^*$ the dual cone in the dual space 
$Im(\partial )^*\otimes\real$.  } 
\end{definition} 

\begin{definition}\label{lambdapi} 
{\rm Let  $\pi \in \Pi$ be any positive path $\pi \in \Pi$.  
We associate to $\pi$ a linear function
$$ \lambda_{\pi} : L(E) \lra \integer$$
by setting $\lambda_{\pi}(e)=1$ if the path $\pi$ crosses the edge $e$, 
and $\lambda_{\pi}(e)=0$ if it doesn't.} 
\end{definition} 

\begin{remark}
{\rm If the path $\pi$ enters a box $b\in B$, then it does 
so by crossing an edge which is part of the opposite corner 
$\cC_b^-$, and it has to leave $b$ by crossing an edge 
which is part of the corner $\cC_b^-$. 
It follows that the corresponding functional $\lambda_{\pi}$ 
is zero on 
$Ker(\partial)=H_1(\Gamma)$, hence it descends to a functional on 
$L(E)/Ker(\partial)=Im(\partial)$, still denoted by $\lambda_{\pi}$. 
By definition,  
$\lambda_{\pi}$ 
is a lattice point in  the dual cone
$C^*\subset Im(\partial )^*\otimes\real$.} 
\end{remark}  
\begin{theorem}\label{cone}
The  semigroup of lattice points in $C^*$ is minimally generated by 
the set of all $\lambda_{\pi}$, where 
$\pi$ runs over the set  $\Pi$ of positive paths.
\end{theorem}
\prf
Let  $\lambda: L(E) \lra \integer $ with $\lambda |_{Ker(\partial)}=0$
and $\lambda(e) \ge 0$ for all $e \in E$. We define the {\em weight} of 
$\lambda$ to be 
\[ w(\lambda) = \sum_{e \in E} \lambda(e). \]
It is clear that $w(\lambda) \geq 0$, and that
$w(\lambda) = 0$ iff $\lambda = 0$. Note also 
that $w(\lambda_{\pi}) = n$ for all $\pi \in \Pi$.

The statement of theorem will be proved if we 
show that $w(\lambda) \geq n$ for 
all non-zero integral linear 
functions $\lambda: L(E) \lra \integer $ with $\lambda |_{Ker(\partial)}=0$
and $\lambda(e) \ge 0$ for all $e \in E$, and, moreover, 
any such $\lambda$  is  a non-negative 
integral linear combination
of $\lambda_{\pi}$ ($\pi \in \Pi$). 

By \ref{basis},  the requirement  $\lambda |_{Ker(\partial)}=0$ 
is equivalent to  $\lambda (\rho_b)=0$ for all $b\in B$, or  
\[  \sum_{e \in \cC_b} \lambda(e) = 
\sum_{e \in \cC_b^-}\lambda(e) \;\; {\rm for\; all\; } b \in B. \]

As in the proof of Proposition \ref{global-sect}, the above condition 
implies that 
if $\lambda \neq 0$, then there exists a roof $\cR_i$  
containing an edge $e$ on which $\lambda$ is nonzero (hence 
$\lambda(e) \geq 1$). 
We start to construct a positive 
path $\pi_i$ from $O_i$ by choosing its edges in such a way 
that $e$ is the first edge of the graph 
$\Gamma$ intersected by $\pi_i$.  
Let $b \in B$ be a box containing $e$ in its opposite corner (i.e., 
$e \in \cC_b^-$). Since  $\lambda (\rho_b)=0$, there must
be an edge $f \in \cC_b$ 
such that $\lambda(f) > 0$ (hence, $\lambda(f) \geq 1$). 
We prolong the path $\pi_i$ through $f$ and enter a next  box $b'$, for which 
$f \in \cC_{b'}^-$. Again there must  exist an edge $g \in \cC_{b'}$ 
such that  $\lambda(g) \geq 1$. $\ldots$ etc. 
Continuing this process, we eventually obtain a positive 
path $\pi_i$, which only
crosses edges $e$ of $\Gamma$ having the property $\lambda(e) \geq 1$. 
This shows that $\lambda':=\lambda-\lambda_{\pi_i}$ is again an integral 
non-negative linear functional on $C$. On the other hand, 
$$w(\lambda') = w(\lambda) - w(\lambda_{\pi_i}) =  w(\lambda) - n.$$
Since $w(\lambda')\geq 0$, this shows that  $w(\lambda) \geq n$. 
By induction on  $w(\lambda)$, we can assume that 
$\lambda'$ is already a non-negative integral linear combination 
of $\lambda_{\pi}$, hence so is $\lambda = \lambda' + \lambda_{\pi_i}$. 
\qed

\vskip 10pt

\begin{definition} 
{\rm We define a partial ordering on the set $\Pi$ of 
positive paths by declaring
$\pi \geq \pi'$ 
if the path $\pi$ runs above the path $\pi'$.} 
\end{definition}
\centerline{Figure 6}

\vskip 20pt
\epsfxsize 8cm
\centerline{\epsfbox{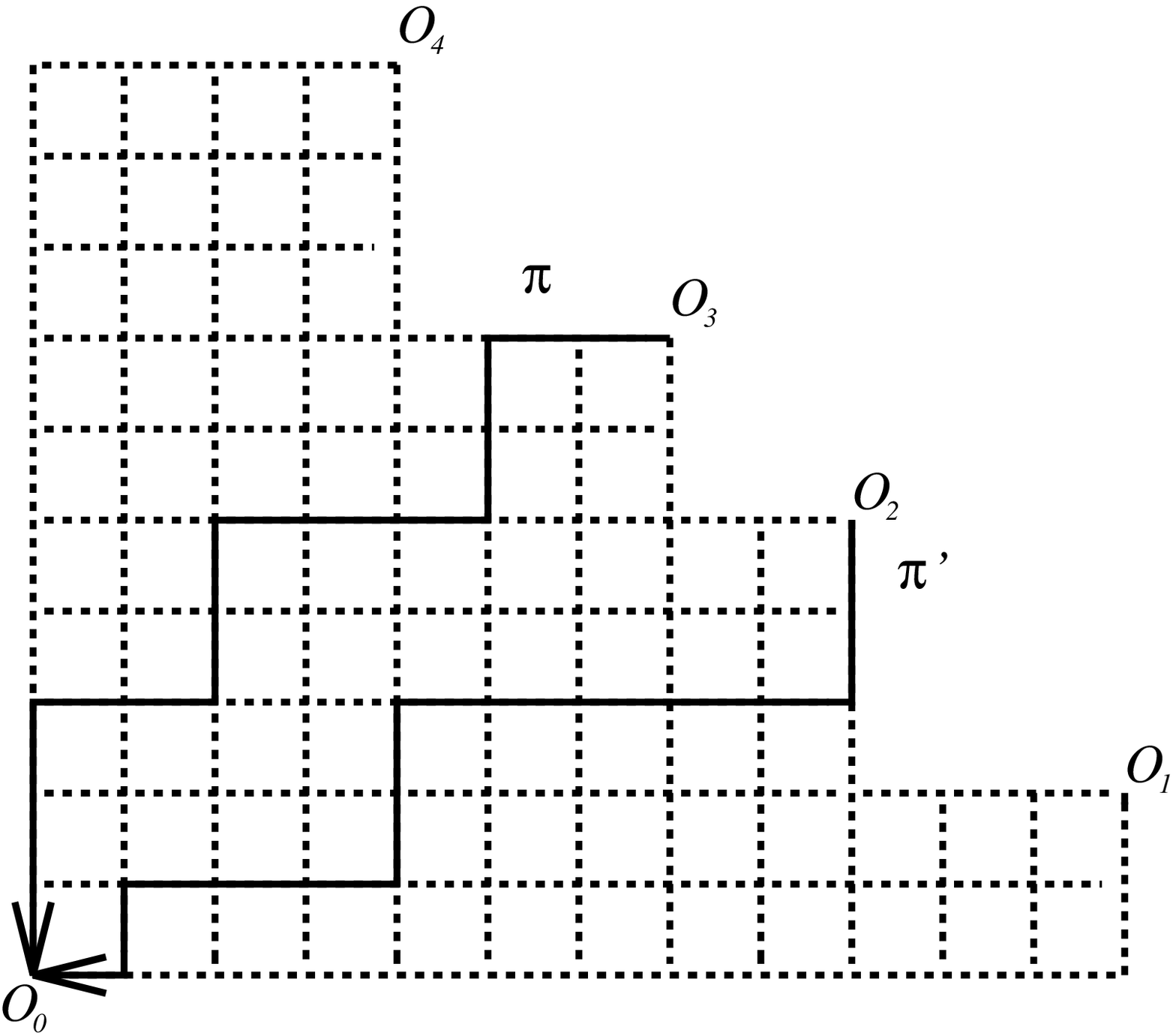}}
\vskip 20pt
\begin{remark}
{\rm It is easy to see that the set $\Pi$ of positive paths together with the above partial ordering is a distributive lattice.  
The maximum
$\max(\pi,\pi')$ for any two paths $\pi$ and $\pi'$ is 
the path bounding the union
of the regions under  $\pi$ and $\pi'$; similarly, $\min(\pi,\pi')$ bounds the
intersection of these regions.
}
\end{remark} 

\begin{definition} 
{\rm Consider the partition of the set of independent variables 
$\{ z_{\pi} \}_{\pi\in \Pi}$ into $l$ disjoint subsets
$$\{z_{\pi} \}_{\pi\in\Pi_i},\; \; i =1, \ldots, l, $$ 
and define $X=X(n_1, \ldots, n_l,n)$ to be the 
subvariety of 
$${\proj}^{N_1-1} \times {\proj}^{N_2-1} \times 
\cdots \times {\proj}^{N_l-1}$$
given by the $l$-homogeneous quadratic equations 
\begin{equation}\label{eq}
z_{\pi}z_{\pi'}-z_{\min(\pi,\pi')}z_{\max(\pi,\pi')} =0,
\end{equation}   
for all pairs of noncomparable elements $\pi,\pi'\in\Pi$}. 
\end{definition}

The variety 
$X$ has been investigated by 
N. Gonciulea and V. Lakshmibai in the papers \cite{GL1,GL2},
where the following result 
has been proved: 

\begin{theorem} {\rm (i)} $X(n_1, \ldots, n_l, n)$  is a $|D|$-dimensional, irreducible, normal, toric variety.

{\rm (ii)} There exists a flat deformation
$$\rho: {\cal X} \lra Spec(\complex[t])$$
such that $\rho^{-1}(0)=X(n_1, \ldots, n_l, n)$, and 
$\rho^{-1}(t)=F(n_1, \ldots, n_l, n) $ for all $t \neq 0$. 
\label{deform}
\end{theorem}

The next theorem describes an  isomorphism
\[  X(n_1, \ldots, n_l, n) \cong P(n_1, \ldots, n_l, n). \]  
  
\begin{theorem} Let  $P= P(n_1, \ldots, n_l, n)$ be the 
toric variety associated with a partial flag manifold $F=F(n_1,\ldots,n_l, n)$.
The line bundles $\cL_i$ $(i =1, \ldots, l)$ 
define an embedding
$$\psi: P   \hookrightarrow {\proj}^{N_1-1} \times 
{\proj}^{N_2-1} \times \cdots \times {\proj}^{N_l-1},$$ 
whose image coincides with the toric variety 
$ X(n_1, \ldots, n_l, n)$.
\label{identif}
\end{theorem}

\prf  We have $X=Proj(\complex[z_{\pi}\; ;\; \pi\in\Pi]/\cI)$, with $\cI$ the ideal generated by the quadratic polynomials in (\ref{eq}), and $Proj$ is taken with respect to the $\integer ^l$-grading given by 
$${\rm deg}(z_{\pi})=\stackrel{i}{(0,\ldots, 0,1,0,\ldots ,0)},\ \ {\rm if}\ \pi\in\Pi_i.$$
If we identify $\Pi$ with the set $\{\; \lambda_{\pi},\pi\in\Pi\;\}\subset Im(\partial)^*$, then $\cI$ is the {\em toric ideal} (see the definition in \cite{sturmfels}, p. 31) associated to this set (this is a standard fact about the ideals associated to distributive lattices; see for example Thm. 4.3 in \cite{GL1} for a proof).
Let $Y$ be the affine toric variety $Spec(\complex[z_{\pi}\; ;\; \pi\in\Pi]/\cI)$. By Theorem \ref{cone} and Prop. 13.5 in \cite{sturmfels},
$Y$ coincides with the affine toric variety defined by the cone $C\subset Im(\partial)\otimes\real$, i.e. $\complex[z_{\pi}\; ;\; \pi\in\Pi]/\cI$ can be identified with the ring $\complex [S_C]$ determined by the semigroup $S_C$ of lattice points in the dual cone $C^*$.

Pick an edge $e_i\in\cR_i$ for each $1\leq i\leq l$, and identify the line bundle $\cL_i$ with $\cO(\sum_{f\in U(e_i)}H_f)$. For each $1\leq i\leq l$, let
$\Delta[e_i]\subset L(D)^*\otimes\real$ be the supporting polyhedron for the global sections of the line bundle $\cL_i$ (cf. Proposition \ref{global-sect});
recall that the lattice points in $\Delta[e_i]$ are given by the linear functions $\pi[e_i]$ ($\pi\in\Pi_i$) defined in the proof of Lemma \ref{cartier}.
Define now for each $i$ a linear function $v[e_i]:L(E)\lra\integer$ by
$$v[e_i](f)=\left\{ \begin{array}{rl}
1,&\mbox {\rm if $f\in U(e)$,}\\
0, &\mbox {\rm otherwise.}\end{array} \right. $$
It is clear that $v[e_i]$ descends to a functional on $Im(\partial)$, and that
for every path $\pi\in\Pi_i$ the functional $\lambda_{\pi}$ in Definition \ref{lambdapi} coincides with $\pi[e_i]+v[e_i]$. For each $1\leq i\leq l$, let $\sigma_i\subset Im(\partial)^*\otimes\real$ be the cone over the translated polyhedron $v[e_i]+\Delta[e_i]$. Then the Minkowski sum 
$\sigma:=\sigma_1+\cdots +\sigma_l$ of these cones coincides with the cone $C^*$, since both $\sigma$ and $C^*$ are generated by the vectors 
$\{ \lambda_{\pi},\; \pi\in\Pi\}$. 
It follows that $P\cong Proj(\complex [S_C])$, where $Proj$ is taken with respect to the natural $\integer ^l$-grading induced by the decomposition
of $C^*$ into the Minkowski sum of the $\sigma_i$'s. 

For each $1\leq i\leq l$, choose an ordering 
$\{\pi_{i,1},\pi_{i,2},\ldots,\pi_{i,N_i}\}$ of the set $\Pi_i$. 
Let $s_{\pi_{i,j}}\in H^0(P,\cL_i)$ denote the section determined by $\pi_{i,j}$ ($i=1,\ldots,l$, $j=1,\ldots,N_i$).
The line bundles $\cL_1,\ldots,\cL_l$ define a morphism
\[\psi :P\lra{\proj}^{N_1-1} \times {\proj}^{N_2-1} \times \cdots \times {\proj}^{N_l-1},\]
\[x\mapsto\left ( [s_{\pi_{1,1}}(x):\ldots :s_{\pi_{1,N_1}}(x)],\ldots,
[s_{\pi_{l,1}}(x):\ldots :s_{\pi_{l,N_l}}(x)]\right ) .\] 
By the above arguments, $\psi$ is the isomorphism 
$$Proj(\complex [S_C])\lra Proj(\complex[z_{\pi}\; ;\; \pi\in\Pi]/\cI), $$
and the theorem is proved.\qed

\vskip 10pt

From \ref{deform} and \ref{identif} we obtain 

\begin{corollary} There exists a flat deformation
$$\rho: {\cal X} \lra Spec(\complex[t])$$
such that $\rho^{-1}(0)=P(n_1, \ldots, n_l, n) 
$ and $\rho^{-1}(t)=F(n_1, \ldots, n_l, n) $ for all $t \neq 0$. 
\end{corollary}

\begin{remark}
{\rm  A description of the singular locus of $P$ was conjectured
by N. Gonciulea and V. Lakshmibai  in the case $F$ a Grassmannian 
(see  \cite{GL2}). Our Theorem \ref{sing} 
proves this  conjecture  and 
its generalization for arbitrary partial flag manifolds $F$.}
\end{remark}

\section{Quantum differential systems} 

\subsection{Quantum $\cD$-module} \label{qdmod}

In order to explain our mirror construction, we give a short overview 
of the quantum cohomology $\cD$-module. The reader is refered to
\cite{G1, K3} for details.

Let $V$ be a smooth projective variety.
Denote $\{ T_a\} _a$ and $\{ T^a\} _a$ two homogeneous bases of
$H^*(V,\rational )$, dual with respect to the Poincar\'e pairing, i.e. 
such that
$$<T_a, T^b>=\delta _{a,b}.$$
We will consider only the even degree part of $H^*(V,\rational )$
and will assume that $H^2(V,\integer )$ and $H_2(V,\integer )$
are torsion-free.
We denote by $1$ the fundamental class of $V$.

To simplify the exposition, suppose there
is a basis $\{ p_i,\ i=1,2,\ldots ,l\}$ of $H^2(V,\integer )$
consisting of nef divisors. 
Let $NE(V)$ be the Mori cone of $V$. 

Introduce formal parameters $q_i$, $i=1,...,l$
and let $\rational [[ q_1,...,q_l]]$ be the ring of formal power series.
The small quantum cohomology ring of 
$V$ will be denoted by $QH^*(V)$. 
This is the free $\rational [[q_1,\ldots,q_l]]$-module
$H^*(V,\rational)\otimes_{\rational}\rational 
[[q_1,\ldots,q_l]]$, together with
a new multiplication given by
$$
T_a\circ T_b=\sum_{\beta \in NE(V)}
\prod_{i=1}^lq_i^{<p_i,\beta >}\left ( \sum _{c}I_{3,\beta}^V(T_aT_bT_c)T^c
\right ) ,
$$
with $I_{3,\beta}^V(T_aT_bT_c)$ the {\em $3$-point, 
genus $0$, Gromov-Witten invariants} of $V$.

\begin{remark} {\rm For the case of a partial flag manifold, the small quantum cohomology ring is well-understood. A presentation of this ring is known
(\cite{AS, K1, K2}), as well as explicit formulas for quantum multiplication (\cite{ICF}).}
\end{remark}

The operators of quantum multiplication with the generators $p_i$ give 
the {\em quantum differential system}, a consistent 
first order partial differential system  (see e.g. \cite{G1}): 

$$
\hbar \frac{\partial}{\partial t_i} \vec{S} = p_i\circ \vec{S},\ i=1,...,l
$$

$$
\hbar \frac{\partial}{\partial t_0} \vec{S} = 1\circ \vec{S},
$$
where $\vec{S}$ is an $H^*(V,\rational )$-valued function in
formal variables $t_0$, and $t_i=\log q_i$, $i=1,...,l$.
Here $\hbar$ is an additional parameter.

Remarkably, a complete set of solutions to this system 
can be written down explicitly in terms of the so-called 
{\em gravitational descendants} \cite{G1}:
$$
\begin{array}{lcl}
\vec{S}_{a}^V:&=&e^{t_0/\hbar}\left ( e^{pt/\hbar}T_a
+\right .\\
 & &\left . +\sum _{\beta \in  NE(V) -0}q^{<p,\beta >}\sum _{b}T^b
\int _{[\overline{M}_{0,2}(V,\beta )]}
\frac{e_1^*(e^{pt/\hbar} T_a)}{\hbar -c}\cup e^*_2(T_b)
\right ) .\\
\end{array}
$$
Here $\overline{M}_{0,2}(V,\beta )$ is 
Kontsevich's space of stable maps, with evaluation morphisms 
$e_1,e_2:\overline{M}_{0,2}(V,\beta )\rightarrow V$ at the two marked points,
$[\overline{M}_{0,2}(V,\beta )]$ is the {\em virtual fundamental class}
(\cite{behrend, li-tian}), and $c$ is the
first Chern class of the line bundle over
$\overline{M}_{0,2}(V,\beta )$ given by the cotangent line
at the first marked point. Finally, $pt$ and 
$q^{<p,\beta>}$ are shorthand notations for $\sum_ip_it_i$ 
and $\prod_iq_i^{<p_i,\beta>}$ respectively.

The {\em quantum $\cD$-module} of $V$ is the $\cD$-module generated by
the functions $<\vec{S},1>$ for all solutions 
$\vec{S}$ to the above differential system.

A general conjecture about the structure of  
quantum $\cD$-modules is Givental's version of the mirror conjecture \cite{G2}:

\begin{conjecture} There exists a family $(M_q,\cF_q,\omega_q)$ 
of (possibly noncompact) complex manifolds $M_q$,
 having the same dimension as $V$, together with 
holomorphic functions $\cF_q$, and
holomorphic volume forms $\omega_q$ such that 
the $\cD$-module generated by integrals
$$\int_{\gamma\subset M_q}e^{(\cF_q+t_0)/\hbar}\omega_q,$$
where $\gamma$ are suitable Morse-theoretic 
middle dimension cycles of the function $Re(\cF_q)$, 
is equivalent to the quantum $\cD$-module of $V$.
\label{mirror}
\end{conjecture}

\subsection{Complete intersections}\label{ci}

Now assume that $V$ is Fano. Let $X$ be the zero 
locus of a generic section of a
decomposable rank $r$ vector bundle 
$$\cE=\bigoplus _{j=1}^rL_j,$$
such that each $L_j$ is generated by global sections. In such a situation one
can also define a quantum ring $QH^*(\cE)$
over the coefficient ring
$\rational [[q_1,...,q_l ]]$, which encodes 
some of the enumerative geometry of rational curves 
on the complete intersection $X$. This leads to a quantum differential 
system for $(V,\cE)$ (see \cite{G3, K3}).
We define degrees of $q_i$'s by requiring that
$$c_1(TV )-c_1(\cE)=\sum(\deg q_i)p_i.$$
Furthermore, we suppose that
all degrees of $q_i$ are nonnegative 
(this is equivalent to the condition that $-K_X$ is nef).
One can write down a similar complete set of 
solutions to the quantum differential system for $(V,\cE)$ \cite{G3}:
$$
\begin{array}{lcl}
\vec{S}_{a}^{\cE}:&=&e^{t_0/\hbar}\left ( e^{pt/\hbar}T_a
+\right .\\
 & &\left . +\sum _{\beta \in  NE(V)  -0}
q^{<p,\beta >}\sum _{b}T^b\int _{[\overline{M}_{0,2}(V,\beta )]}
\frac{e_1^*(e^{pt/\hbar} T_a)}{\hbar -c}\cup e^*_2(T_b)
\cup E_{\beta }\right ) ,\\
\end{array}
$$
where $E_{\beta}$ is the Euler class of 
the vector bundle on $\overline{M}_{0,2}(V,\beta )$ 
whose fibre over a point $(C,\mu ; x_1,x_2)$ is 
the subspace of $H^0(\mu^*\cE)$ consisting of sections 
vanishing at $x_2$, and the rest of notations are as above.  
 
Consider the cohomology valued functions 
$$S_V:=\sum _a <\vec{S}_{a}^V ,1> T^a$$
and
$$S_{\cE}:=\sum _a <\vec{S}_{a}^{\cE} ,c_r(\cE)> T^a.$$
These functions are given explicitly by the expressions
$$ 
S_V= e^{(t_0+pt)/\hbar}\left(
1+\sum _{\beta \in  NE(V) - 0}
q^{<p,\beta >}(e_1)_*
\left (\frac{1}{\hbar -c}\right )\right ) 
$$
and
$$ 
S_{\cE}= e^{(t_0+pt)/\hbar}\left(
c_r(\cE)+\sum _{\beta \in  NE(V) - 0}
q^{<p,\beta >}(e_1)_*
\left (\frac{e_1^*(c_r(\cE))\cdot E'_{\beta}}{\hbar -c}\right )\right ) ,
$$
where now $E'_{\beta}$ is the Euler class of 
the vector bundle on $\overline{M}_{0,2}(V,\beta )$ 
whose fibre over a point $(C,\mu ; x_1,x_2)$ is 
the subspace of $H^0(\mu^*\cE)$ consisting of sections 
vanishing at $x_1$.

\begin{remark}{\rm If we view $X$ as an abstract variety, the general theory in \S\ref{qdmod} gives a $H^*(X,\rational)$-valued function $S_X$. The functions $S_X$ and $S_{\cE}$ are closely related. For example, if 
$i^*:H^2(V,\integer )\tilde{\rightarrow} H^2(X,\integer)$, where $i:X\hookrightarrow V$
is the inclusion, then $i_*(S_X)=S_{\cE}$.}
\end{remark}
Now consider a new cohomology valued function

$$
I_{\cE}=e^{(t_0+pt)/\hbar}\left ( c_r(\cE)+
\sum _{\beta \in  NE(V)  -0}q^{<p,\beta >}
\prod _j\prod _{m=0}^{<c_1(L_j),\beta >}
(c_1(L_j) + m\hbar ) (e_1)_*
(\frac{1}{\hbar -c})\right).
$$
In general it is very hard to compute $S_X$ or $S_{\cE}$ explicitly. However, note that $I_{\cE}$ can be computed directly from the function $S_V$ associated to the ambient manifold, which in many cases turns out to be more tractable.
It is therefore extremely useful to have a result relating $S_{\cE}$ and $I_{\cE}$.
Extending ideas of Givental, B. Kim \cite{K3} has recently proved the following theorem,
which applies to the cases considered in this paper:    
\begin{theorem} \label{qlef}
If $V$ is a homogeneous space and $X\subset V$ is the zero locus of a generic section of a nonnegative decomposable vector bundle $\cE$, 
then $S_{\cE}$ and $I_{\cE}$ coincide up to a weighted homogeneous triangular 
change of variables:
$$t_0\lra t_0+f_0(q)\hbar+f_{-1}(q),\ \ \  \log q_i\lra \log q_i+f_i(q),
\;\;\; i=1,\dots ,l,$$
where $f_{-1},f_0,f_1,\ldots ,f_l$ are weighted 
homogeneous formal power series supported in $NE(V)  -0$, 
with $\deg f_{-1}=1$, and $\deg f_i=0$, $i=0,1,\dots ,l$.
\end{theorem}

In particular, this implies that the coefficient $\Phi_V$ of the cohomology class $1\in H^*(V,\rational)$ 
in $S_V$, and the coefficient $\Phi_X$ of $c_r(\cE)$
in $I_{\cE}$ (specialized to $\hbar=1$, $t_0=0$) are related in a very simple way. Namely, if 
$$\Phi_V =  \sum_{\beta \in  NE(V) -0} a_{\beta} q^{<p,\beta >}$$
$$\Phi_X = \sum_{\beta \in  NE(V) -0} b_{\beta} q^{<p,\beta >},$$ 
then 
\begin{equation}\label{trick}
b_{\beta }=a_{\beta}\prod_{i=1}^r (<c_1(L_i), \beta >!).
\end{equation}

We will refer to Theorem \ref{qlef} as the {\em quantum hyperplane section theorem}. The relation $(\ref{trick})$ above was called the ``trick with factorials'' in \cite{grass}.

\section{The mirror construction}

In this section we give a partially conjectural mirror construction for partial
flag manifolds, and use it to obtain an explicit hypergeometric series as the 
power series expansion of the integral representation. The case of Calabi-Yau
complete intersections is then discussed in some detail. 
 
\subsection{Hypergeometric solutions for partial flag manifolds}

Let $F= F(n_1, \ldots, n_l, n)$ be a partial flag manifold. In the  
notations of Section \ref{toric}, we   
introduce  $l$ independent variables 
$q_i$, $i=1,2,\ldots,l$ (each $q_i$ corresponds to the  roof $\cR_i$), 
$|B|$ independent variables $\tilde{q}_b$,  $b \in B$, and 
$|E|$ independent variables $y_e$, $ e \in E$. 
Consider the following set of algebraically independent polynomial 
equations:\\
$\bullet$ Roof equations: For $i=1,2,\ldots, l$
\begin{equation}\label{rfeqn} 
\cF_i:=\prod_{e \in \cR_i}y_e-q_i = 0.
\end{equation}
$\bullet$ Box equations: For $b=\{e,f,g,h\} \in B$ 
\begin{equation}\label{boxeqn}
\cG_b:=y_e y_f -\tilde{q}_by_g y_h=0,
\end{equation}
where $\{e,f\} = \cC_b$.

This set of equations was discussed by Givental \cite{G2}, and was used to give an integral representation for
the solutions to the quantum cohomology differential equations for the special case of complete flag manifolds. The results in that paper were the starting point for our investigations. We describe below Givental's result and our (conjectural) generalization to a general partial flag manifold.

Let $ {\bf A}^{\mid E \mid}$ be the complex affine space with the coordinates 
$y_e$ $( e \in E)$. 
For fixed parameters values of 
$$(q, \tilde{q}) :=(q_1,\ldots,q_l,\ldots,\tilde{q}_b,\ldots)$$
we obtain an affine  variety: 
$$M_{q,\tilde{q}}:=\{u \in {\bf A}^{\mid E \mid}\;\;|\;\;\cF_i=0,\; i=1,\ldots ,l,\; \cG_b=0,\; b\in B\}.$$
If all components
of $(q, \tilde{q})$ are non-zero, $M_{q, \tilde{q}}$ is isomorphic to the 
torus $({\bf C}^*)^{\mid D \mid}$.

One can define on $M_{q, \tilde{q}}$ a holomorphic volume form
$$\omega_{q, \tilde{q}} :=Res_{M_{q, \tilde{q}}}\left 
(\frac{\Omega}{\prod_{i=1}^l\cF_i\prod_{b \in B}\cG_b}
\right ) ,$$ 
where $$\Omega:=\wedge_{e \in E}dy_e.$$
Let $\cF=\sum_{e \in E}y_e$.
Consider the integral
$$I_{\gamma}(q, \tilde{q}):=\int_{\gamma}e^{\cF} \omega_{q, \tilde{q}},$$
where $\gamma \in H_{|D|} ({M_{q, \tilde{q}}}, Re(\cF)=-\infty)$.   
We put
$$\Phi_{\gamma}(q_1,\ldots,q_l):=I_{\gamma}(q_1,\ldots,q_l, 1,1,\ldots,1).$$

We can now formulate a precise version of conjecture \ref{mirror}:

\begin{conjecture} Let $\vec{S}$ be any solution to 
the quantum differential system for $F$. 
Then the component $<\vec{S},1>$ can be expressed as $\Phi_{\gamma}(q)$
for some $\gamma \subset M_{q,1}$.
\label{flagmirror}
\end{conjecture}

\begin{remark}\label{complete}
{\rm This conjecture generalizes Givental's mirror theorem for complete flag manifolds \cite{G2}. } 
\end{remark} 

\begin{definition}
{\rm Let $W$ denote the set of edges in the diagram $\Lambda$ that intersect $\Gamma$. We orient the vertical edges in $W$ upwards and the horizontal edges to the right. Let $V:=B\cup\{ 0,1,2,\ldots ,l\}$.
For $w\in W$, the {\em tail} $t(w)$ of $w$, is defined 
to be the box $b_1\in B$ where $w$ starts. Similarly, the 
{\em head} $h(w)$ of $w$ is the box $b_2\in B$ where $w$ ends. 
If $w$ crosses the roof $\cR_i$, so that its \lq\lq head" 
is outside the graph $\Gamma$, we put $h(w):=i$, and if the \lq\lq tail"
of $w$ is outside $\Gamma$, we put $t(w)=0$.
In the sense of duality of planar graphs, the graph with vertices $V$, edges
$W$, and incidence given by $h,t:W\lra V$ is {\em dual} to the graph $\Gamma$ with all stars collapsed to one point.} 
\end{definition}

\begin{definition}
{\rm For each cone $\sigma\in\widehat{\Sigma}$ of maximal dimension we define a cycle $\gamma=\gamma_{q,\tilde {q}}(\sigma)$ in $M_{q,\tilde{q}}$ by
$$\gamma :=\{\; \mid y_e\mid =1,\; {\rm for\; all}\; e\in E\; {\rm with}\; \delta(e)\in\sigma\} .$$ 
Note that the $y_f$'s with $\delta(f)\notin\sigma$ are determined uniquely
by the $y_e$'s with $\delta(e)\in\sigma$ and the roof and box equations $(\ref{rfeqn})$, $(\ref{boxeqn})$.}
\end{definition}
 
The cycle $\gamma$ is a real torus, of dimension equal to ${\rm dim}_{\complex}(M_{q,\tilde q})={\rm dim}_{\complex}(F)$. Since it is defined over the entire family of $M_{q,\tilde{q}}$'s, it is invariant under monodromy.
The integral over this special cycle will be denoted by $I(q,\tilde{q})$. 

\begin{definition}
{\rm The specialization $\Phi_F(q):=I(q_1,\dots,q_l,1,\ldots,1)$ is called
the {\em hypergeometric series of the partial flag manifold $F$.}}
\end{definition}

It turns out that $I(q,\tilde{q})$ has a nice 
power series expansion.

\begin{theorem}\label{hyperg}
$I(q, \tilde{q})=\displaystyle{\sum_{m_1,\ldots,m_l,\ldots,m_b,\ldots}}
A_{m_1,\ldots ,m_l,\ldots,m_b,\ldots}\; 
q_1^{m_1}\ldots q_l^{m_l}\prod_{b \in B}{\tilde{q}_b}^{m_b}$,
with 
\[ A_{m_1,\ldots,m_l,\ldots,m_b,\ldots}:=
\frac{1}{(m_1!)^{k_1+k_2}}\frac{1}{(m_2!)^{k_2+k_3}}
\ldots \frac{1}{(m_l!)^{k_{l}+k_{l+1}}}B_{m_1,\ldots,m_l,\ldots,m_b,\ldots},\]
 
\[ B_{m_1,\ldots,m_l,\ldots,m_b,\ldots}:=\prod_{w \in W}{m_{h(w)} 
\choose m_{t(w)}}. \] 

\end{theorem} 
\prf By Leray's theorem, the integral is equal to 
$$\int_{T(\gamma_{q,\tilde{q}}(\sigma))}e^{\cF}\frac{\Omega}{\prod_{i=1}^l\cF_i \prod_{b \in B}\cG_b},$$
where $T$ is the tube map. For $\mid q\mid <1$, $\mid\tilde{q}\mid <1$, the cycle $T(\gamma_{q,\tilde{q}}(\sigma))$ is homologous to the cycle 
$$
T:=\{ y\in {\bf A}^{\mid E \mid}\; |\; \mid y_e\mid =1,\; {\rm for\; all}\; e\in E\} 
$$
in the complement of the hypersurfaces $y_e=0$.
We now expand all the terms in the integrand:
$$
e^{\cF}=\sum_{d=0}^{\infty}\frac{1}{d!}{\cF}^d=\sum_{d_e \ge
0}\frac{\prod_{e \in E}y_e^{d_e}}{\prod_{e \in E} d_e !}
$$
$$
\frac{1}{\cF_i}=\frac{1}{\prod_{e \in \cR_i}y_e}\sum_{m_i\ge 0}
\left (\frac{q_i}{\prod_{e \in \cR_i}y_e}\right )^{m_i}
$$
$$
\frac{1}{\cG_b}=\frac{1}{y_ey_f}\sum_{m_b \ge 0}\left
(\frac{\tilde{q_b}y_gy_h}{y_ey_f}\right )^{m_b},
$$
where $\{e,f\}$ makes up the corner and $\{g,h\}$ the opposite corner of the
box $b=\{e,f,g,h\}$.
The integral picks up precisely the constant coefficient of the
following power
series in the $y_e$'s, with parameters the $q$'s and $\tilde{q}$'s:
$$
\sum_{d_e, m_i, m_b \ge 0} \frac{\prod_{e \in E}y_e^{d_e}}{\prod_{e \in E}
d_e !}\prod_{i=1}^l\left (\frac{q_i}{\prod_{e \in \cR_i}y_e}\right )^{m_i}
\prod_{b \in B}\left (\frac{\tilde{q_b}y_gy_h}{y_ey_f}\right )^{m_b}
$$
Now there are three types of edges.
\vskip 10pt
Type I: $e \in \cR_i$ for some $i=1,\ldots,l$. Then $e$ also is edge of
the opposite corner of a unique box $b$. Only the terms with
$$ d_e=m_i-m_b$$ will give a contribution.
\vskip 10pt
Type II: $e \in b \cap b'$, for two boxes $b$ and $b'$. We can then
assume that $e$ is part of the corner of $b$, and the opposite corner of
$b'$. Only the terms with 
$$ d_e=m_b-m_{b'}$$
will give a contribution.
\vskip 10pt
Type III: $e$ is contained in a unique $b \in B$. In this case
$e$ is part of the corner of $b$. Only the terms with
$$ d_e=m_b$$
will give a contribution. 
\vskip 10pt

Hence we see that the integral is given by
the series:
$$ 
\sum_{m_i, m_b \ge 0}\frac{1}{\prod_{e \in
E}d_e!}\prod_{i=1}^lq_i^{m_i}\prod_{b \in B}\tilde{q_b}^{m_b},
$$
where for each edge the number $d_e$ is determined by the $m_i$ and $m_b$
by the above equations.
We can rewrite this coefficient nicely in terms of binomial coefficients as
follows. Each edge $w \in W$ of the diagram $\Lambda$ intersects
precisely one edge $e \in E$ of type $I$ or type $II$. The corresponding coefficient $d_e$ is then given by
$$
d_e=m_{h(w)}-m_{t(w)}
$$ 
Trivially,
$$
\frac{1}{\prod_{e \in E}d_e!}=\frac{\prod_{w \in W}m_{h(w)}!}{\prod_{w \in W}m_{h(w)}!\prod_{e \in E}d_e!}
$$
The heads of arrows $w \in W$ which are {\em not} tails are the heads of arrows
intersecting the edges of type I. The tails of arrows $w \in W$ which
are not heads are in bijection to the edges of type III. Hence, when we
pull out a factor $\prod_{i=1}^l\prod_{e \in \cR_i}(m_e!)$ from the
denominator of the left hand side of the above equality, the other terms
in numerator and denominator can precisely be combined into the product
$$\prod_{w \in W}{m_{h(w)}\choose m_{t(w)}}.$$
This proves the result. \qed

\begin{remark} 
{\rm Note that $I(q, \tilde{q})$ is the generalized 
hypergeometric series for the smooth toric variety
$\widehat{P}$ defined in Section 3. The parameters $q_1,\ldots ,q_l$ correspond to the generators of $Pic (\widehat{P})$ coming from the singular variety $P$
(the pull-backs of the line bundles $\cL_1,\ldots ,\cL_l$),
while $\tilde {q}_b$ correspond to the additional generators of  
$Pic (\widehat{P})$.}
\end{remark}

Theorem \ref{hyperg} shows that it is very easy to write down the power series expansion for $I(q,\tilde{q})$ directly from the diagram.

\begin{example} 
{\rm 
$F(2,5)$ (the Grassmannian of $2$-planes in $\complex^5$)
\vskip 40pt
\epsfxsize 5cm
\centerline{\epsfbox{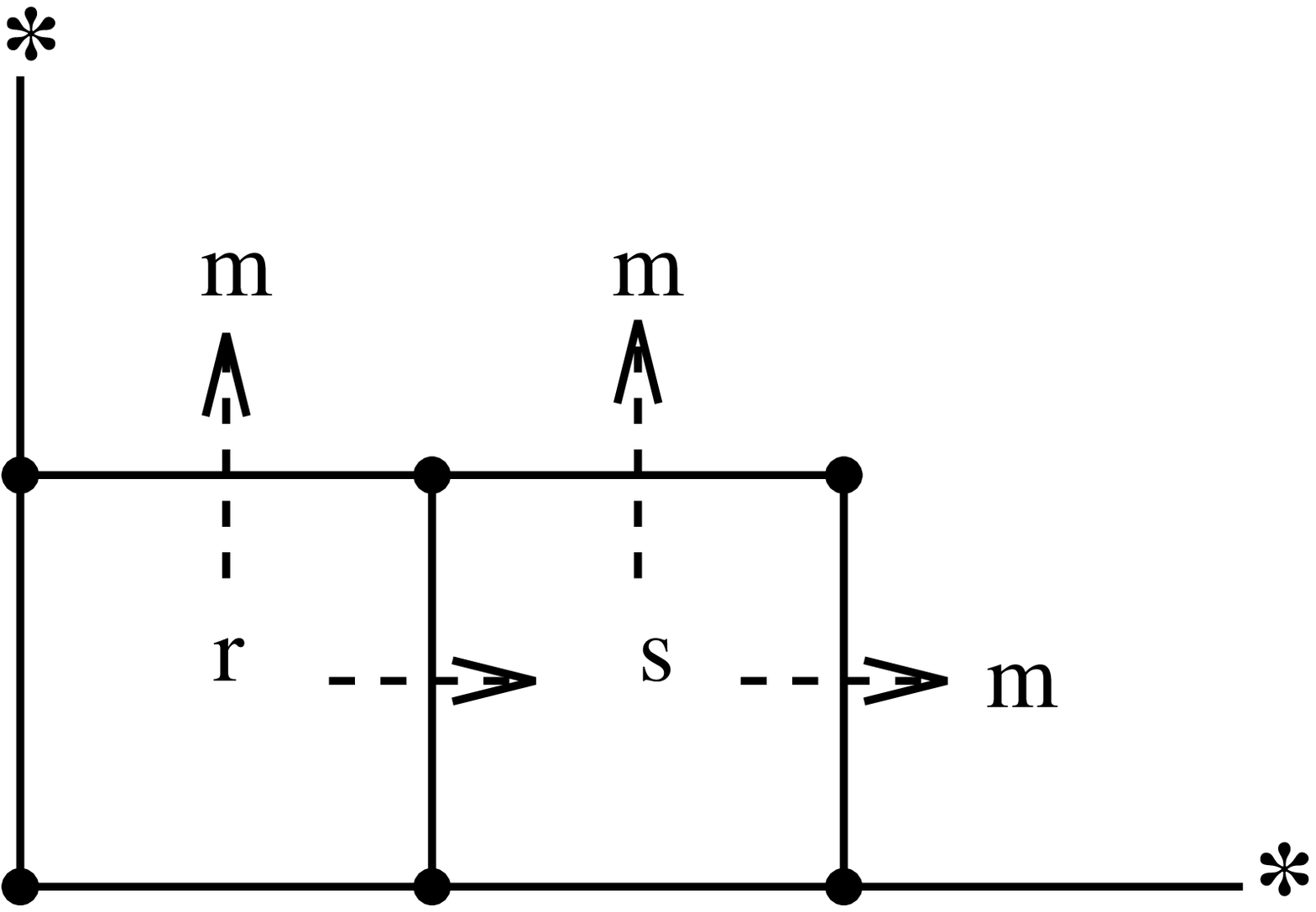}}
\vskip 30pt
Hence we read off:
$$I(q, \tilde{q})=\sum_{m,r,s \ge 0}\frac{1}{(m!)^5}{s \choose r}{m \choose r}
{m \choose s}^2q^m \tilde{q}_1^r\tilde{q}_2^s.$$
} 
\end{example}


\begin{example}
{\rm 
$F(3,6)$ (the Grassmannian of $3$-planes in $\complex^6$)
\vskip 20pt
\epsfxsize 5cm
\centerline{\epsfbox{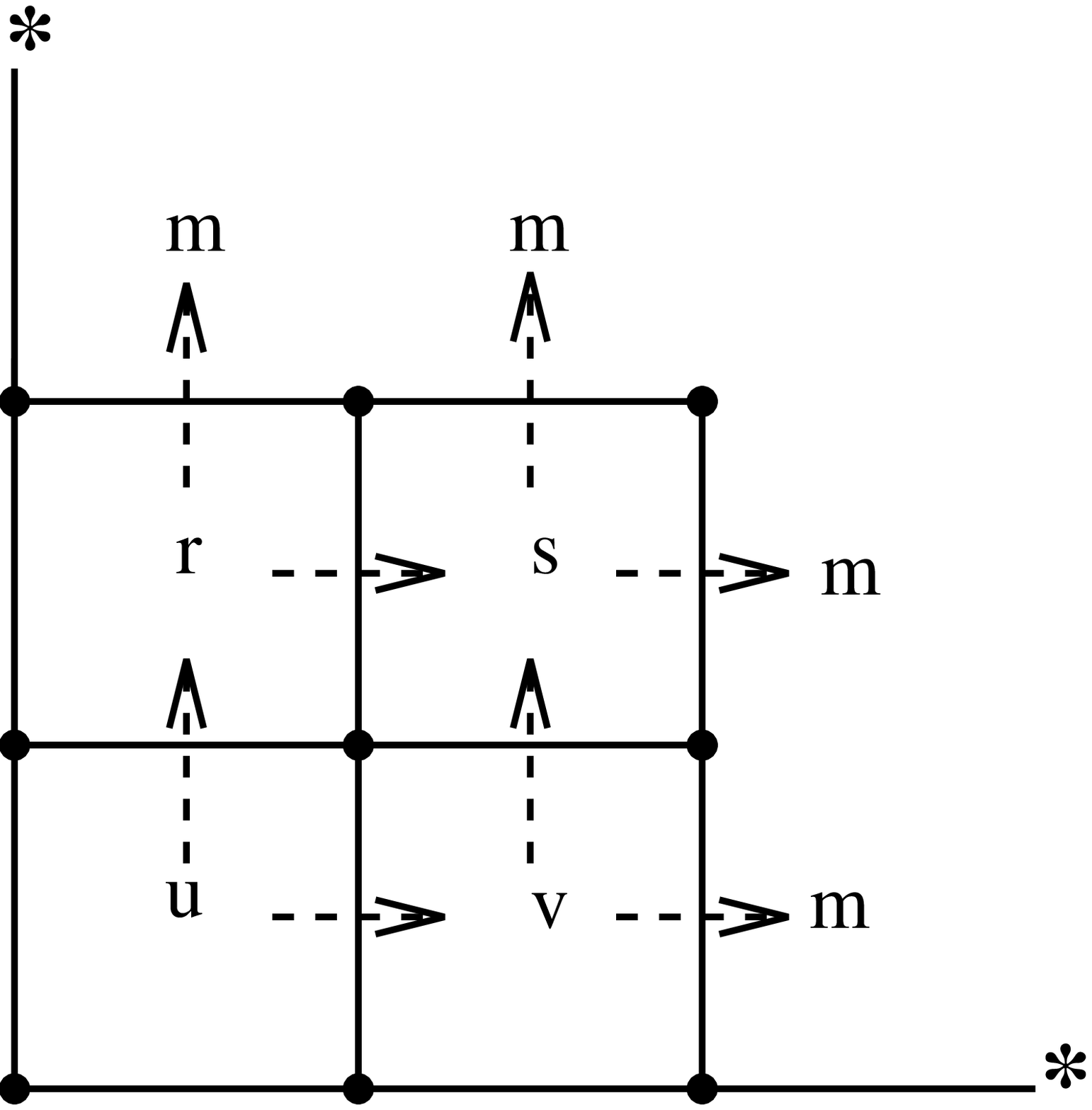}}
\vskip 20pt
Hence we read off:

$$I(q, \tilde{q})=\sum_{m,r,s,u,v}\frac{1}{(m!)^6}{r \choose u}
{v \choose u}{s \choose r}{s \choose v}{m \choose r}{m\choose s}^2 
{m \choose v} q^m \tilde{q}_1^r \tilde{q}_2^s \tilde{q}_3^u \tilde{q}_4^v.$$
}
\end{example}

\begin{example} 
{\rm 
$F(1,2,3,4)$ (the variety of complete flags in $\complex ^4$)
\vskip 20pt
\epsfxsize 5cm
\centerline{\epsfbox{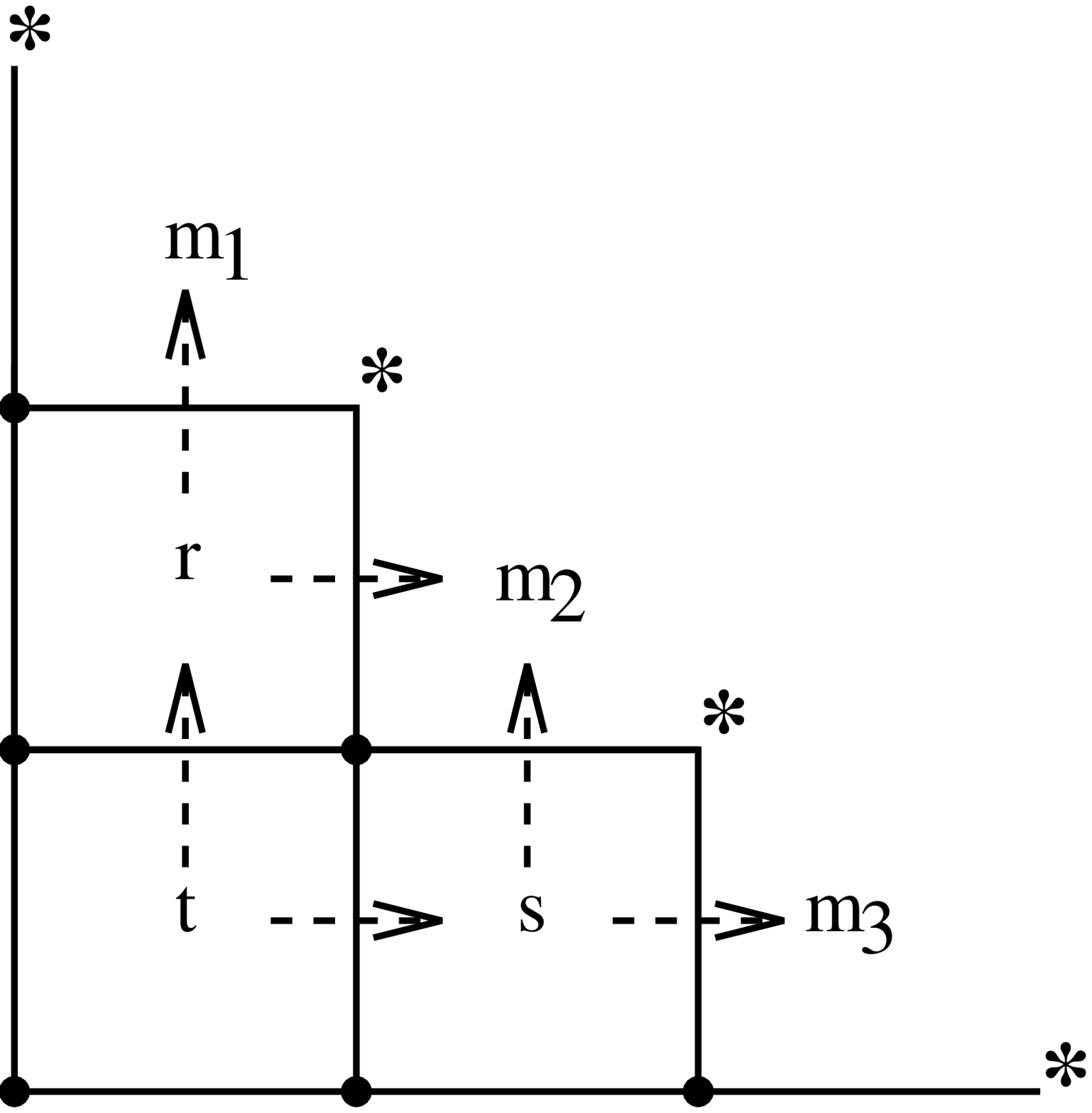}}
\vskip 20pt

Hence we read off:
$$I(q,\tilde{q})=\sum_{m_1,m_2,m_3,r,s,t}A_{m_1,m_2,m_3,r,s,t}q_1^{m_1} 
q_2^{m_2} q_3^{m_3} \tilde{q}_1^r \tilde{q}_2^s \tilde{q}_3^t, $$
with
$$A_{m_1,m_2,m_3,r,s,t}=\frac{1}{(m_1!)^2(m_2!)^2(m_3!)^2}{r \choose t}
{s \choose t}{m_1 \choose r}{m_2 \choose r} {m_2 \choose s}{m_3 \choose s}.$$ 
}
\end{example}
\vskip 20pt

A weaker version of Conjecture \ref{flagmirror} is

\begin{conjecture}\label{S0}
The series $\Phi_F:=I(q,1)$ is the coefficient of the cohomology class 1 in the $H^*(F,\rational)$-valued function $S_F$ describing the quantum $\cD$-module of $F$, that is,
$$\Phi_F=1+\sum_{\overline m:=(m_1,\ldots ,m_l)\neq 0}\left (\int_{\overline{M}_{0,2}(F,\overline m)}
\frac{e_1^*(e^{Ct} \Omega_F)}{1-c}\cup e_2^*(1)\right )q_1^{m_1}\ldots q_l^{m_l},$$ 
where $Ct$ stands for $C_1t_1+\cdots +C_lt_l$, with $\{ C_1,\ldots ,C_l\}$ the Schubert basis of $H^2(F,\rational)$, and $\Omega_F$ is
the cohomolgy class of a point.
\end{conjecture}

\begin{remark} $(i)$ {\rm Besides the case of complete flag manifolds (cf. Remark \ref{complete}), there is another case for which the above conjecture agrees with previously known results. Consider the partial flag manifold
$F:=F(1,n-1,n)$ of flags $V^1\subset V^{n-1}\subset\complex ^n$. The Pl\" ucker
embedding, identifies $F$ with a $(1,1)$ hypersurface in 
$\proj ^{n-1}\times\proj ^{n-1}$. The hypergeometric series for  
$\proj ^{n-1}\times\proj ^{n-1}$ is
$$\sum_{m_1 ,m_2\geq 0}\frac{1}{(m_1!)^n(m_2!)^n}q_1^{m_1}q_2^{m_2},$$
(cf. \cite{G3}), and by the quantum hyperplane section theorem (\cite{G3}, \cite{K3}) we obtain that the hypergeometric series for $F$ is
\begin{equation}\label{known}
\sum_{m_1 ,m_2\geq 0}\frac{(m_1+m_2)!}{(m_1!)^n(m_2!)^n}q_1^{m_1}q_2^{m_2}.
\end{equation}
On the other hand, the recipe of Theorem \ref{hyperg} gives
the formula
\begin{equation}\label{new}
\sum_{m_1 ,m_2\geq 0}\frac{\sum_s{m_1\choose s}{m_2\choose s}}{(m_1!)^{n-1}(m_2!)^{n-1}}q_1^{m_1}q_2^{m_2}
\end{equation}
for the hypergeometric series of $F$. 
The identity $\sum_s{m_1\choose s}{m_2\choose s}={m_1+m_2\choose m_1}$ implies  that the series $(\ref{known})$ and $(\ref{new})$ coincide.
}

$(ii)$ {\rm The quantum Pieri formula \cite{ICF} gives explicitly the quantum product of a special Schubert class with a general one, and in particular the quantum product of a Schubert divisor with any other Schubert class. Using this, one can write down in reasonably low dimensional cases the 
quantum differential system for $F$ and reduce this first order system to higher order differential equations satisfied by the components. 
In particular, one can write down the differential operators annihilating the
component $<\vec{S},1>$ of any solution $\vec{S}$, and check by direct computation that the hypergeometric series $\Phi_F(q)$ of Theorem \ref{hyperg} is annihilated by these operators. In \cite{grass} this is done for the Grassmannians containing complete intersection Calabi-Yau $3$-folds.
For the complete flag manifolds, the operators are known to be the operators for the quantum Toda lattice (see \cite{Kim4})}
\end{remark}

\subsection{Calabi-Yau complete intersections in $F(n_1,\ldots,n_l,n)$}

Recall that ${\rm Pic}(F)$ is generated by the line bundles $\cO (C_i)$, $i=1,\ldots ,l$, which also generate the (closed) K\"ahler cone. Hence any line bundle $\cH$ on $F$ which is globally generated is of the form $\cO(\overline{d}):=\cO (\sum_{i=1}^l d^{(i)}C_i)$, with $d^{(i)}$ nonnegative. 
The common zero locus of $r$ general sections of line bundles $\cO(\overline{d}_1),
\ldots ,\cO(\overline{d}_r)$ will be denoted by
$X:=X_{\overline {d}_1,\ldots ,\overline {d}_r}$. 

Assuming Conjecture \ref{S0}, it follows from the quantum hyperplane section theorem that the hypergeometric series $\Phi_X$
has the expression
\begin{equation}\label{cyhyperg}
\Phi_X=\sum_m\prod_{j=1}^r(\sum_{i=1}^ld_j^{(i)}m_i)!A_{m_1,\ldots ,m_l}q_1^{m_1}\ldots q_l^{m_l},
\end{equation}
where $A_{m_1,\ldots ,m_l}$ are the coefficients of $\Phi_F$ in Theorem \ref{hyperg}.

From now on the complete intersection $X_{\overline {d}_1,\ldots ,\overline {d}_r}$ is assumed to be a Calabi-Yau manifold. The construction of mirrors  described in \cite{grass} for the case when $F$ is a Grassmannian can be extended to the case of a general $F$ as follows:

$X$ can be regarded as the intersection of $F\subset \proj^{N_1-1}\times\ldots\times\proj^{N_l-1}$ with $r$ general hypersurfaces
$Z_j\; (j=1,\ldots ,r)$ in $\proj^{N_1-1}\times\ldots\times\proj^{N_l-1}$, with 
$Z_j$ of multidegree $(d_j^{(1)},\ldots ,d_j^{(l)})$. Let $Y$ be the Calabi-Yau complete intersection of the same hypersurfaces with the toric degeneration $P$
of $F$.

For each edge $e\in \cup_{i=1}^l\cR_i$ , which is part of a roof, define polynomials 
$$\varphi_e(y):=\sum_{f\in U(e)}c_fy_f,$$
where $c_f$ are generically chosen complex numbers
(recall that we have defined $U(e)$ as the set consisting of $e$, 
together with all
edges in the graph $\Gamma$ which are either directly below $e$, if $e$ is horizontal, or directly to the left of $e$, if $e$ is vertical).

\vskip 10pt

Partition each of the roofs $\cR_i$, $i=1,\ldots ,l$ into 
$r$ disjoint subsets
$$\cR_i=\cR_{i,1}\cup\ldots\cup\cR_{i,r}$$
such that $\mid \cR_{i,j}\mid=d_i^{(j)}$. It follows from Definition \ref{linebundles} and Theorem \ref{identif} that the toric Weil divisor
$$\sum_{i=1}^l\sum_{e\in\cR_{i,j}}H_e$$ is Cartier, and
$$\cO(\sum_{i=1}^l\sum_{e\in\cR_{i,j}}H_e)\cong \cL_1^{\otimes d_j^{(1)}}\otimes\ldots\otimes\cL_l^{\otimes d_j^{(l)}}.$$

\vskip 10pt

Consider the torus $T$ in the affine space 
$\cong{\bf A}^{\mid E \mid}$ given by the following set of equations:

$\bullet$ Roof equations: For $i=1,2,\ldots, l$
$$ \prod_{e \in \cR_i}y_e=1.$$ 

$\bullet$ Box equations: For $b=\{e,f,g,h\} \in B$ 
$$y_e y_f -y_g y_h=0,$$
where $\{e,f\}$ form the corner $\cC_b$ of $b$.

Introduce additional independent variables $x_d,\; d\in D$, one for each generator of the lattice $L(D)$. For every edge $e\in E$, set 
$$x^{\delta(e)}:=x_{h(e)}(x_{t(e)})^{-1},$$
where, as before, $h(e)$ (resp. $t(e)$) is the head (resp. tail) of $e$.
The torus $T$ can be identified with $Spec(\complex[x_d,x_d^{-1}\; ;\; d\in D])$, with the embedding $T\hookrightarrow{\bf A}^{\mid E\mid}$ induced by the ring homomorphism
$$\complex[y_e\; ;\; e\in E]\longrightarrow\complex[x_d,x_d^{-1}\; ;\; d\in D],\;\;\; y_e\mapsto x^{\delta(e)}.$$
With this identification, we obtain Laurent polynomials
$$\varphi_e(x):=\sum_{f\in U(e)}c_fx^{\delta(f)}.$$
For $j=1,\ldots,r$, let $\nabla_j$ be the Newton polyhedron of the Laurent polynomial
$$\cP_j:=1-\sum_{i=1}^l\sum_{e\in\cR_{i,j}}\varphi_e(x).$$
The polyhedra $\nabla_j,\; j=1,\ldots,r$ define a nef-partition of the anticanonical class of ${P}$ (see definitions in \cite{LB,BB2}), and according to \cite{BS} and \cite{LB}, the mirror family $Y^*$ of the Calabi-Yau 
complete intersection $Y \subset{P}$ consists of Calabi-Yau compactifications of the general complete intersections in $T$ defined by the equations
\begin{equation}\label{mirroreq}
1-\sum_{i=1}^l\sum_{e\in\cR_{i,j}}\varphi_e(x)=0,\;\;\; j=1,\dots,r.
\end{equation}

\begin{conjecture}\label{CY} Let $Y_0^*$ be a Calabi-Yau compactification
of a general complete intersection in $T$ defined by the equations $(\ref{mirroreq})$, 
with the additional requirement that the coefficients satisfy the relation $$c_{f_1}c_{f_2}=c_{f_3}c_{f_4}$$
whenever $\{ f_1,f_2,f_3,f_4\}$ make up a box $b\in B$, with $\{ f_1,f_2\}$
forming the corner $\cC_b$ of $b$. Then a minimal desingularization of $Y_0^*$
is a mirror of a generic complete intersection Calabi-Yau $X\subset F$.
\end{conjecture}

The main period of the mirror $Y^*$ of $Y$ is given by
$$ \Phi_{Y}=\int_{\gamma}Res_{M_q,\tilde{q}}\left (\frac{\Omega}{\prod_{j=1}^r\cE_j\prod_{i=1}^l\cF_i\prod_{b \in B}\cG_b}
\right ) ,$$ 
where the extra factors $\cE_j$ come from the nef-partition
of the anti-canonical class of ${P}$ described above. 
Specifically,$$\cE_j:=1-\sum_{i=1}^l\sum_{e\in\cR_{i,j}}\sum_{f\in U(e)}y_f,\ \ \ j=1,\ldots,r.$$ 
By direct expansion of the integral defining $\Phi_Y$ (as in Theorem \ref{hyperg}), followed by the specialization $\tilde{q}_b=1$, $b\in B$, one gets exactly the hypergeometric series $\Phi_X$.

\vskip 10pt

Finally, we discuss some applications to the case when $X\subset F$ is a Calabi-Yau $3$-fold. 

First, as discussed in \cite{grass}, our construction can be interpreted via {\em conifold transitions}. Indeed, by Theorem \ref{sing}, if $X$ is generic, 
then its degeneration $Y\subset P$
is a singular Calabi-Yau $3$-fold, whose singular locus consists of finitely many nodes. The resolution of singularities $\widehat{P}\lra P$ induces a small resolution $\widehat{Y}\lra Y$. In other words the (nonsingular) Calabi-Yau's $X$ and $\widehat Y$ are related by a conifold transition, and Conjecture \ref{CY} essentially states that their mirrors are related in a similar fashion.

Second, it is well understood (see e.g. \cite{BS}) that the knowledge of the hypergeometric series $\Phi_X$ for a Calabi-Yau $3$-fold gives the virtual numbers of rational curves on $X$ via a formal calculation. In \cite{grass}
we have used the hypergeometric series $(\ref{cyhyperg})$ to compute these numbers for complete intersections in Grassmannians.

\subsection{List of Calabi-Yau complete intersection $3$-folds}
\vskip 10pt
Recall that if $F:=F(n_1,\dots ,n_{l},n)$ is a partial flag manifold, then 
\begin{equation}\label{dimension} 
{\rm dim}\, (F)= \sum_{i=1}^{l}(n_i-n_{i-1})(n-n_i). 
\end{equation} 
In the Schubert basis of the Picard group, the anticanonical bundle of $F$ 
is given by
\begin{equation}\label{-K} 
\omega_F^{-1}=\cO \left (\sum_{i=1}^l (n_{i+1}-n_{i-1})C_i\right ). \end{equation}

A (general) complete intersection Calabi-Yau 3-fold in $F$ is 
the common zero locus of $r:={\rm dim}\, (F)-3$ general sections  $s_j\in H^0(F,\cO(\overline{d}_j))$, where $\cO(\overline{d}_j),\, j=1,2,\ldots ,r$ are line bundles with $\bigotimes_{j=1}^r\cO(\overline{d}_j)=\omega_F^{-1}$.
Hence, if $F$ contains a complete intersection Calabi-Yau $3$-fold, then necessarily
\begin{equation}\label{constraint} 
{\rm dim}\, (F)\leq 3+\sum_{i=1}^l (n_{i+1}-n_{i-1})=n+n_l-n_1+3. 
\end{equation}

\begin{proposition}\label{bound}
If $F:=F(n_1,\dots ,n_{l},n)$ is a partial flag manifold
containing a complete intersection Calabi-Yau $3$-fold, and $F$ is not a projective space, or one of the manifolds $F(1,n-1,n)$, then $n\leq 7$.
\end{proposition}

\prf Using (\ref{dimension}), after some manipulation, one can rewrite the inequality (\ref{constraint}) as
\begin{equation}\label{cons2}
(n_1-1)(n-n_1-1)+(n_2-n_1)(n-n_2-1)+\ldots+(n_l-n_{l-1})(n-n_l-1)\leq 4.
\end{equation}
There are two cases.

\noindent {\it Case 1}: $n_1>1$. Then it is easy to see that
$(n_1-1)(n-n_1-1)>4$ for $n\geq 8$, unless $n_1=n-1$, in which case $F$
is a projective space.

\noindent {\it Case 2}: $n_1=1$. If $l=1$, then $F$ is a projective space, 
so we may assume $l\geq 2$. As above, 
$(n_2-1)(n-n_2-1)>4$ for $n\geq 8$, unless $n_2=n-1$, in which case $F=F(1,n-1,n)$.\qed

\begin{remark} {\rm The flag manifold $F(1,n-1,n)$ sits as a $(1,1)$ hypersurface in
$\proj^{n-1}\times\proj^{n-1}$. Hence these cases (as well as the case when $F$ is projective space) can be viewed as particular instances of complete intersection Calabi-Yau's in toric varieties.}  
\end{remark}

We list below all the partial flag manifolds (not excluded by proposition \ref{bound}) for which the inequality (\ref{constraint}) is satisfied. 
The anticanonical class of $F$, denoted by $-K_F$, is expressed in terms of the natural Schubert basis of the Picard group. The last column of the table below
contains the possible splittings of the anticanonical class into 
${\rm dim}\, (F)-3$ nonnegative divisors.

In general, there is a natural duality isomorphism 
\begin{equation}\label{iso}
F(n_1,\ldots ,n_l,n)\cong F(n-n_l,\ldots ,n-n_1,n).
\end{equation}
This is taken into account by listing only one of the two isomorphic flag manifolds.
It may also be that the flag manifold is self-dual, i.e. (\ref{iso}) is an
automorphism, and two families of
complete intersection Calabi-Yau $3$-folds corresponding to different splittings of the anticanonical class are interchanged by the duality automorphism.
Whenever this happens (e.g., when $F$ parametrizes complete flags), only one of the two splittings of $-K_F$ is listed.

{\tiny \begin{center}

\begin{tabular}{|c|c|c|c|c|} \hline
 & & & & \\
$n$ & $F$ & ${\rm dim}\, (F)$ & $-K_F$ & splitting of $-K_F$ \\
 & & & & \\
\hline  & & & & \\
$7$ & $F(2,7)$ & 10 & $7$ & $7(1)$\\  & & & & \\
\hline  & & & & \\
$7$  & $F(1,2,7)$ & 11 & $(2,6)$ & $2(1,0)+6(0,1)$\\  & & & & \\
\hline  & & & & \\
$7$  & $F(1,5,7)$ & 14 & $(5,6)$ & $5(1,0)+6(0,1)$\\  & & & & \\
\hline  & & & & \\
$7$  & $F(1,2,6,7)$ & 15 & $(2,5,5)$ & $2(1,0,0)+5(0,1,0)+5(0,0,1)$\\ & & & & \\
\hline  & & & & \\
$6$ & $F(2,6)$ & 8 & $6$ & $(2)+4(1)$\\ & & & & \\
\hline  & & & & \\
$6$  & $F(3,6)$ & 9 & $6$ & $6(1)$\\ & & & & \\
\hline & & & & \\
$6$  & $F(1,2,6)$ & 9 & $(2,5)$ & $(2,0)+5(0,1)$\\
 & & & & $(1,0)+(1,1)+4(0,1)$ \\
 & & & & $2(1,0)+(0,2)+3(0,1)$ \\ & & & & \\
\hline  & & & & \\
$6$ & $F(1,3,6)$ & 11 & $(3,5)$ & $3(1,0)+5(0,1)$ \\ & & & & \\
\hline & & & & \\
$6$  & $F(1,4,6)$ & 11 & $(4,5)$ & $(2,0)+2(1,0)+5(0,1)$ \\ 
 & & & & $3(1,0)+(1,1)+4(0,1)$ \\
 & & & & $4(1,0)+(0,2)+3(0,1)$ \\ & & & & \\
\hline  & & & & \\
$6$  & $F(1,2,5,6)$ & 12 & $(2,4,4)$ & $(2,0,0)+4(0,1,0)+4(0,0,1)$ \\
 & & & & $(1,0,0)+(1,1,0)+3(0,1,0)+4(0,0,1)$ \\
 & & & & $(1,0,0)+(1,0,1)+4(0,1,0)+3(0,0,1)$ \\
 & & & & $2(1,0,0)+(0,2,0)+2(0,1,0)+4(0,0,1)$ \\
 & & & & $2(1,0,0)+(0,1,1)+3(0,1,0)+3(0,0,1)$ \\
 & & & & $2(1,0,0)+4(0,1,0)+(0,0,2)+2(0,0,1)$ \\ & & & & \\
\hline & & & & \\
$6$  & $F(1,3,5,6)$ & 13 & $(3,4,3)$ & $3(1,0,0)+4(0,1,0)+3(0,0,1)$\\ & & & & \\ \hline & & & & \\
$5$ & $F(2,5)$ & 6 & $5$ & (3)+2(1) \\
 & & & & $2(2)+(1)$ \\ & & & & \\
\hline & & & & \\
$5$  & $F(1,2,5)$ & 7 & $(2,4)$ &  $(2,0)+(0,2)+2(0,1)$  \\
 & & & & $(1,0)+(1,1)+(0,2)+(0,1)$ \\
 & & & & $2(1,1)+2(0,1)$ \\
 & & & & $2(1,0)+2(0,2)$ \\
 & & & & $(1,0)+(1,2)+2(0,1)$ \\
 & & & & $(2,1)+3(0,1)$ \\ & & & & \\
\hline & & & & \\
$5$  & $F(2,3,5)$ & 8 & $(3,3)$ & $(1,0)+(2,0)+3(0,1)$ \\
 & & & & $2(1,0)+(1,1)+2(0,1)$ \\ & & & & \\
\hline & & & & \\
$5$  & $F(1,3,5)$ & 8 & $(3,4)$  & $(3,0)+4(0,1)$ \\
 & & & & $(1,0)+(2,1)+3(0,1)$ \\
 & & & & $(1,1)+(2,0)+3(0,1)$ \\
 & & & & $(1,0)+2(1,1)+2(0,1)$ \\
 & & & & $(1,0)+(2,0)+(0,2)+2(0,1)$ \\
 & & & & $3(1,0)+2(0,2)$ \\
 & & & & $3(1,0)+(0,1)+(0,3)$ \\
 & & & & $2(1,0)+(1,2)+2(0,1)$ \\ & & & & \\
\hline

\end{tabular}
\end{center}

\begin{center}

\begin{tabular}{|c|c|c|c|c|} \hline
 & & & & \\
$n$ & $F$ & ${\rm dim}\, (F)$ & $-K_F$ & splitting of $-K_F$ \\
 & & & & \\
\hline  & & & & \\
$5$  & $F(1,2,4,5)$ & 9 & $(2,3,3)$ & $2(1,0,0)+(0,3,0)+3(0,0,1)$ \\
 & & & & $2(1,0,0)+3(0,1,0)+(0,0,3)$ \\
 & & & & $(2,1,0)+2(0,1,0)+3(0,0,1)$ \\
 & & & & $(2,0,1)+3(0,1,0)+2(0,0,1)$ \\
 & & & & $(1,2,0)+(1,0,0)+(0,1,0)+3(0,0,1)$ \\
 & & & & $2(1,0,0)+(0,2,1)+(0,1,0)+2(0,0,1)$ \\
 & & & & $2(1,0,0)+(0,1,2)+2(0,1,0)+(0,0,1)$ \\
 & & & & $(1,0,0)+(1,0,2)+3(0,1,0)+(0,0,1)$ \\
 & & & & $(1,1,1)+(1,0,0)+2(0,1,0)+2(0,0,1)$ \\
 & & & & $(2,0,0)+(0,2,0)+(0,1,0)+3(0,0,1)$ \\
 & & & & $(2,0,0)+(0,0,2)+3(0,1,0)+(0,0,1)$ \\
 & & & & $2(1,0,0)+(0,2,0)+(0,1,0)+(0,0,2)+(0,0,1)$ \\
 & & & & $2(1,1,0)+(0,1,0)+3(0,0,1)$ \\
 & & & & $(1,1,0)+(1,0,1)+2(0,1,0)+2(0,0,1)$ \\
 & & & & $(1,1,0)+(1,0,0)+(0,1,1)+(0,1,0)+2(0,0,1)$ \\
 & & & & $2(1,0,0)+2(0,1,1)+(0,1,0)+(0,0,1)$ \\
 & & & & $(1,0,0)+(1,0,1)+(0,1,1)+2(0,1,0)+(0,0,1)$ \\
 & & & & $2(1,0,1)+3(0,1,0)+(0,0,1)$ \\
 & & & & $(2,0,0)+(0,1,1)+2(0,1,0)+2(0,0,1)$ \\
 & & & & $(1,1,0)+(0,2,0)+(1,0,0)+3(0,0,1)$ \\
 & & & & $(1,0,1)+(0,2,0)+(1,0,0)+(0,1,0)+2(0,0,1)$ \\
 & & & & $2(1,0,0)+(0,2,0)+(0,1,1)+2(0,0,1)$ \\
 & & & & $(1,1,0)+(1,0,0)+2(0,1,0)+(0,0,2)+(0,0,1)$ \\
 & & & & $(1,0,1)+(1,0,0)+3(0,1,0)+(0,0,2)$ \\
 & & & & $2(1,0,0)+(0,1,1)+2(0,1,0)+(0,0,2)$ \\ & & & & \\
\hline & & & & \\
$5$  & $F(1,2,3,5)$ & 9 & $(2,2,3)$ & $(2,0,0)+2(0,1,0)+3(0,0,1)$ \\
 & & & & $(0,2,0)+2(1,0,0)+3(0,0,1)$ \\
 & & & & $(0,0,2)+2(1,0,0)+2(0,1,0)+(0,0,1)$ \\
 & & & & $(1,1,0)+(1,0,0)+(0,1,0)+3(0,0,1)$ \\
 & & & & $(1,0,1)+(1,0,0)+2(0,1,0)+2(0,0,1)$ \\
 & & & & $(0,1,1)+2(1,0,0)+(0,1,0)+2(0,0,1)$ \\ & & & & \\
\hline & & & & \\
$5$  & $F(1,2,3,4,5)$ & 10 & $(2,2,2,2)$ & $(2,0,0,0)+2(0,1,0,0)+2(0,0,1,0)+2(0,0,0,1)$ \\
 & & & & $2(1,0,0,0)+(0,2,0,0)+2(0,0,1,0)+2(0,0,0,1)$ \\
 & & & & $(1,1,0,0)+(1,0,0,0)+(0,1,0,0)+2(0,0,1,0)+2(0,0,0,1)$ \\
 & & & & $(1,0,0,1)+(1,0,0,0)+2(0,1,0,0)+2(0,0,1,0)+(0,0,0,1)$ \\
 & & & & $(1,0,1,0)+(1,0,0,0)+2(0,1,0,0)+(0,0,1,0)+2(0,0,0,1)$ \\
 & & & & $(0,1,1,0)+2(1,0,0,0)+(0,1,0,0)+(0,0,1,0)+2(0,0,0,1)$ \\ & & & & \\
\hline  & & & & \\
$4$ & $F(2,4)$ & 4 & $4$  & $(4)$ \\ & & & & \\
\hline & & & & \\
$4$  & $F(1,2,4)$ & 5 & $(2,3)$  & $(1,0)+(1,3)$\\
 & & & & $ (1,1)+(1,2)$ \\
 & & & & $ (2,1)+(0,2)$ \\
 & & & & $ (2,2)+(0,1)$ \\ & & & & \\
\hline & & & & \\
$4$  & $F(1,2,3,4)$ & 6 & $(2,2,2)$ & $(2,0,0)+(0,2,0)+(0,0,2)$ \\
 & & & & $(1,1,0)+(1,0,1)+(0,1,1)$ \\
 & & & & $(1,2,0)+(1,0,0)+(0,0,2)$ \\
 & & & & $(1,2,0)+(1,0,1)+(0,0,1)$ \\
 & & & & $(2,1,0)+(0,1,0)+(0,0,2)$ \\
 & & & & $(2,1,0)+(0,1,1)+(0,0,1)$ \\
 & & & & $(2,0,1)+(0,2,0)+(0,0,1)$ \\
 & & & & $(2,0,1)+(0,1,1)+(0,1,0)$ \\
 & & & & $2(1,1,0)+(0,0,2)$ \\
 & & & & $2(1,0,1)+(0,2,0)$ \\
 & & & & $(2,2,0)+2(0,0,1)$ \\
 & & & & $(2,0,2)+2(0,1,0)$ \\ & & & & \\
\hline 

\end{tabular}
\end{center}
}


\end{document}